# ASYMPTOTIC ANALYSIS OF MULTISCALE APPROXIMATIONS TO REACTION NETWORKS

By Karen Ball,[1] Thomas G. Kurtz,[2] Lea Popovic and Greg Rempala[3]

*IDA Center for Communications Research, University of Wisconsin, University of Minnesota and University of Louisville*

A reaction network is a chemical system involving multiple reactions and chemical species. Stochastic models of such networks treat the system as a continuous time Markov chain on the number of molecules of each species with reactions as possible transitions of the chain. In many cases of biological interest some of the chemical species in the network are present in much greater abundance than others and reaction rate constants can vary over several orders of magnitude. We consider approaches to approximation of such models that take the multiscale nature of the system into account. Our primary example is a model of a cell's viral infection for which we apply a combination of averaging and law of large number arguments to show that the "slow" component of the model can be approximated by a deterministic equation and to characterize the asymptotic distribution of the "fast" components. The main goal is to illustrate techniques that can be used to reduce the dimensionality of much more complex models.

**1. Stochastic models for reaction networks.** A reaction network is a chemical system involving multiple reactions and chemical species. The simplest stochastic model for a network treats the system as a continuous time Markov chain whose state $X$ is a vector giving the number of molecules of each species present with each reaction modeled as a possible transition for the state. The model for the $k$th reaction is determined by a vector of inputs $\nu_k$ specifying the number of molecules of each chemical species that are

Received July 2005; revised February 2006.
[1]Supported in part by VIGRE Grant, University of Indiana.
[2]Supported in part by NSF Grant DMS-05-03983.
[3]Supported in part by CGeMM Intramural Grant, University of Louisville.

consumed in the reaction, a vector of outputs $\nu'_k$ specifying the number of molecules of each species that are created in the reaction, and a function of the state $\lambda_k(x)$ that gives the rate at which the reaction occurs. Specifically, if the reaction occurs at time $t$, the new state becomes

$$X(t) = X(t-) + \nu'_k - \nu_k.$$

Let $R_k(t)$ denote the number of times that the $k$th reaction occurs by time $t$. Then the state of the system at time $t$ can be written as

$$X(t) = X(0) + \sum_k R_k(t)(\nu'_k - \nu_k) = (\nu' - \nu)R(t),$$

where $\nu'$ is the matrix with columns given by the $\nu'_k$, $\nu$ is the matrix with columns given by the $\nu_k$, and $R(t)$ is the vector with components $R_k(t)$.

$R_k$ is a counting process with intensity $\lambda_k(X(t))$ (called the *propensity* in the chemical literature) and can be written as

$$R_k(t) = Y_k\left(\int_0^t \lambda_k(X(s))\,ds\right),$$

where the $Y_k$ are independent unit Poisson processes. Note that writing $R_k$ in this form makes it clear why $\lambda_k$ is referred to as a rate.

Defining $|\nu_k| = \sum_i \nu_{ik}$, the stochastic form of the law of mass action says that the rate should be given by

$$\lambda_k^N(x) = \kappa_k \frac{\prod_i \nu_{ik}!}{N^{|\nu_k|-1}} \binom{x}{\nu_{1k} \cdots \nu_{mk}} = N\kappa_k \frac{\prod_i \nu_{ik}!}{N^{|\nu_k|}} \binom{x}{\nu_{1k} \cdots \nu_{mk}},$$

where $N$ is a scaling parameter usually taken to be the volume of the system times Avogadro's number, and $\kappa_k$ is a constant specifying the rate of the reaction. Note that the rate is proportional to the number of distinct subsets of the molecules present that can form the inputs for the reaction. Intuitively, this assumption reflects the idea that the system is *well stirred*, in the sense that all molecules are equally likely to be at any location at any time.

1.1. *Law of large numbers and diffusion approximations.* If $N$ is the volume times Avogadro's number and $x$ gives the number of molecules of each species present, then $c = N^{-1}x$ gives the concentrations in moles per unit volume. With this scaling and a large volume (where large can be pretty small since Avogadro's number is $6 \times 10^{23}$),

$$\lambda_k^N(x) \approx N\kappa_k \prod_i c_i^{\nu_{ik}} \equiv N\tilde{\lambda}_k(c). \tag{1.1}$$

Since the law of large numbers for the Poisson process implies $N^{-1}Y(Nu) \approx u$, (1.1) implies

$$C(t) = N^{-1}X(t) \approx C(0) + \sum_k \int_0^t \kappa_k \prod_i C(s)_i^{\nu_{ik}} (\nu'_k - \nu_k)\,ds,$$

which in the large volume limit gives the classical deterministic law of mass action

$$\dot{C}(t) = \sum_k \kappa_k \prod_i C(t)_i^{\nu_{ik}} (\nu_k' - \nu_k).$$

Similarly, since an appropriately renormalized Poisson process can be approximated by a standard Brownian motion, that is,

$$\frac{Y(Nu) - Nu}{\sqrt{N}} \approx W(u),$$

we can derive a diffusion approximation for the Markov chain by replacing $Y_k(Nu)$ by $\sqrt{N}W_k(u) + Nu$, that is,

$$\begin{aligned} C^N(t) &= C^N(0) + \sum_k N^{-1} Y_k\biggl(\int_0^t \lambda_k(X^N(s))\,ds\biggr)(\nu_k' - \nu_k) \\ &\approx C^N(0) + \sum_k N^{-1/2} W_k\biggl(\int_0^t \tilde{\lambda}_k(C^N(s))\,ds\biggr)(\nu_k' - \nu_k) \\ &\quad + \int_0^t F(C^N(s))\,ds, \end{aligned}$$

where

$$F(c) = \sum_k \tilde{\lambda}_k(c)(\nu_k' - \nu_k).$$

The diffusion approximation is given by the equation

$$\begin{aligned} \tilde{C}^N(t) &= \tilde{C}^N(0) + \sum_k N^{-1/2} W_k\biggl(\int_0^t \tilde{\lambda}_k(\tilde{C}^N(s))\,ds\biggr)(\nu_k' - \nu_k) \\ &\quad + \int_0^t F(\tilde{C}^N(s))\,ds, \end{aligned}$$

which is distributionally equivalent to the Itô equation

$$\begin{aligned} \tilde{C}^N(t) &= \tilde{C}^N(0) + \sum_k N^{-1/2} \int_0^t \sqrt{\tilde{\lambda}_k(\tilde{C}^N(s))}\,d\tilde{W}_k(s)(\nu_k' - \nu_k) \\ &\quad + \int_0^t F(\tilde{C}^N(s))\,ds \\ &= \tilde{C}^N(0) + \sum_k N^{-1/2} \int_0^t \sigma(\tilde{C}^N(s))\,d\tilde{W} \\ &\quad + \int_0^t F(\tilde{C}^N(s))\,ds, \end{aligned}$$

where $\sigma(c)$ is the matrix with columns $\sqrt{\tilde{\lambda}_k(c)}(\nu_k' - \nu_k)$. A precise version of this approximation is given in [7]. (See also [4], Chapter 10, [5], Chapter 7, and [12].)

1.2. *Multiscale approximations.* Interest in modeling chemical reactions within cells has led to renewed interest in stochastic models, since the number of molecules involved, at least for some of the species, may be sufficiently small that the deterministic model does not provide a good representation of the behavior of the system. Modeling is further complicated by the fact that some species may be present in much greater abundance than others. In addition, the rate constants $\kappa_k$ may vary over several orders of magnitude. With these two issues in mind, we consider a different approach to deriving a scaling limit approximation of the model.

$N$ will still denote a scaling parameter for the model, but it is no longer interpreted in terms of volume or Avogadro's number. In fact, $N^{-1}$ plays the same role as $\varepsilon$ in a perturbation analysis of a deterministic model (see, e.g., [11]). $N$ may have no physical meaning, but it will have a specific (hopefully large) value in any physical or biological setting in which the approximation is applied.

For example, let $N$ be of the order of magnitude of the abundance of the most abundant species in the system. For each species $i$, we then specify a parameter $0 \le \alpha_i \le 1$ and normalize the number of molecules by $N^{\alpha_i}$, defining

$$Z_i(t) = N^{-\alpha_i} X_i(t).$$

$\alpha_i$ should be selected so that $Z_i = O(1)$, but that still leaves a degree of arbitrariness regarding the selection. Note that $\alpha_i$ could be zero, so $Z_i$ could still be integer-valued.

We want to express the reaction rates in terms of $Z$ rather than $X$ and also to take into account large variation in the reaction rates. Consequently, we introduce another set of exponents $\beta_k$ for the reactions and now assume that the reaction rates can be written as $N^{\beta_k}\lambda_k(z)$, where $\lambda_k(z) = O(1)$ for all relevant values of $z$. The model becomes

$$Z_i(t) = Z_i(0) + \sum_k N^{-\alpha_i} Y_k\left(\int_0^t N^{\beta_k}\lambda_k(Z(s))\,ds\right)(\nu_k' - \nu_k).$$

Our goal is to derive simplified models under the assumption that $N$ is large, where "large" may be much smaller than Avogadro's number. We demonstrate that this process may lead to interesting and reasonable models by analyzing a number of examples in the literature.

1.3. *Outline of the paper.* Reaction networks of interest in biology can be very high dimensional involving many chemical species and many reactions. Consequently, there has been considerable effort to exploit the multiscale nature of these systems to derive reduced models. In Section 2 we borrow examples from a number of these papers to illustrate how the kind of scaling limits we have in mind can be used to provide a rigorous and intuitive approach to model reduction. The primary focus of the paper is a model of an intracellular viral infection given in [10] and studied further in [6]. We analyze this model in Section 3, and give a systematic identification of the scaling parameters in Section A.2.

## 2. Examples.

2.1. *Simple crystallization.* We consider a model studied by Haseltine and Rawlings [6] using the parameters in Table I of their paper. The system involves four species and two reactions:

$$2A \xrightarrow{\kappa_1} B, \qquad A + C \xrightarrow{\kappa_2} D.$$

The model satisfies

$$\begin{aligned} X_A(t) &= X_A(0) - 2Y_1\left(\int_0^t \tfrac{1}{2}\kappa_1 X_A(s)(X_A(s)-1)\,ds\right) \\ &\quad - Y_2\left(\int_0^t \kappa_2 X_A(s) X_C(s)\,ds\right), \\ X_B(t) &= X_B(0) + Y_1\left(\int_0^t \tfrac{1}{2}\kappa_1 X_A(s)(X_A(s)-1)\,ds\right), \\ X_C(t) &= X_C(0) - Y_2\left(\int_0^t \kappa_2 X_A(s) X_C(s)\,ds\right). \end{aligned}$$

Following Rawlings and Haseltine, $X_A(0) = 10^6$, $X_B(0) = 0$, $X_C(0) = 10$, and $\kappa_1 = \kappa_2 = 10^{-7}$. Let $N = 10^6$, and take $\alpha_A = \alpha_B = 1$ and $\alpha_C = 0$. Writing $\kappa_1 = \kappa_2 = 0.1 \times N^{-1}$, the normalized system becomes

$$\begin{aligned} Z_A^N(t) &= 1 - N^{-1} 2Y_1\left(N\int_0^t 0.05 Z_A^N(s)(Z_A^N(s) - N^{-1})\,ds\right) \\ &\quad - N^{-1} Y_2\left(\int_0^t 0.1 Z_A^N(s) Z_C^N(s)\,ds\right), \\ Z_B^N(t) &= N^{-1} Y_1\left(N\int_0^t 0.05 Z_A^N(s)(Z_A^N(s) - N^{-1})\,ds\right), \\ Z_C^N(t) &= 10 - Y_2\left(\int_0^t 0.1 Z_A^N(s) Z_C^N(s)\,ds\right). \end{aligned}$$

Letting $N \to \infty$, the simplified system is

$$Z_A(t) = 1 - \int_0^t 0.1 Z_A(s)^2\,ds,$$

$$Z_B(t) = \int_0^t 0.05 Z_A(s)^2\,ds,$$

$$Z_C(t) = 10 - Y_2\biggl(\int_0^t 0.1 Z_A(s) Z_C(s)\,ds\biggr),$$

which gives

$$Z_A(t) = \frac{1}{1 + 0.1t}.$$

Since $Z_A$ is deterministic, $Z_C$ is a linear death process with time-varying rate $\lambda(t) = 0.1 Z_A(t)$. Consequently, for any $t > 0$, the distribution of $Z_C(t)$ is Binomial$(10, p(t))$ with

$$p(t) = \exp\biggl\{-\int_0^t 0.1 Z_A(s)\,ds\biggr\} = \frac{1}{1 + 0.1t}.$$

In particular,

$$E[Z_C(t)] = \frac{10}{1 + 0.1t}, \qquad \mathrm{Var}[Z_C(t)] = \frac{t}{(1 + 0.1t)^2},$$

compare favorably with the simulation results in Figure 1 of [6].

2.2. *Enzyme kinetics.* Rao and Arkin [9] analyze a model of enzyme kinetics, involving an enzyme, substrate, their enzyme-substrate complex and a product of this complex

$$E + S \underset{\kappa_{-1}}{\overset{\kappa_1}{\rightleftarrows}} ES, \qquad ES \xrightarrow{\kappa_2} P + E.$$

The state of this system can be represented by

$$\begin{aligned}
X_{\mathrm{s}}(t) &= X_{\mathrm{s}}(0) + Y_{-1}\biggl(\int_0^t \kappa_{-1} X_{\mathrm{es}}(r)\,dr\biggr) \\
&\qquad - Y_1\biggl(\int_0^t \kappa_1 X_{\mathrm{e}}(r) X_{\mathrm{s}}(r)\,dr\biggr), \\
X_{\mathrm{es}}(t) &= X_{\mathrm{es}}(0) - Y_{-1}\biggl(\int_0^t \kappa_{-1} X_{\mathrm{es}}(r)\,dr\biggr) \\
&\qquad + Y_1\biggl(\int_0^t \kappa_1 X_{\mathrm{e}}(r) X_{\mathrm{s}}(r)\,dr\biggr) - Y_2\biggl(\int_0^t \kappa_2 X_{\mathrm{es}}(r)\,dr\biggr), \\
X_{\mathrm{e}}(t) &= X_{\mathrm{e}}(0) + Y_{-1}\biggl(\int_0^t \kappa_{-1} X_{\mathrm{es}}(r)\,dr\biggr)
\end{aligned}$$

$$- Y_1\left(\int_0^t \kappa_1 X_{\mathrm{e}}(r)X_{\mathrm{s}}(r)\,dr\right) + Y_2\left(\int_0^t \kappa_2 X_{\mathrm{es}}(r)\,dr\right),$$

$$X_{\mathrm{p}}(t) = Y_2\left(\int_0^t \kappa_2 X_{\mathrm{es}}(r)\,dr\right),$$

where $X_{\mathrm{s}}$ gives the number of substrate molecules, $X_{\mathrm{e}}$ the number of enzymes, $X_{\mathrm{es}}$ the number of enzyme complexes, and $X_{\mathrm{p}}$ the number of molecules of the reaction product. Following Rao and Arkin, take $X_{\mathrm{s}}(0) = 100$, $X_{\mathrm{e}}(0) = 1000$, $\kappa_1 = \kappa_{-1} = 1$, $\kappa_2 = 0.1$. Let $N = 1000$, and define $Z_{\mathrm{s}}^N = N^{-2/3}X_{\mathrm{s}}$, $Z_{\mathrm{es}}^N = N^{-2/3}X_{\mathrm{es}}$, $Z_{\mathrm{p}}^N = N^{-2/3}X_{\mathrm{p}}$, $Z_{\mathrm{e}}^N = N^{-1}X_{\mathrm{e}}$, and $\kappa_2 = N^{-1/3}$. Then the normalized system becomes

$$\begin{aligned}
Z_{\mathrm{s}}^N(t) &= 1 + N^{-2/3}Y_{-1}\left(\int_0^t N^{2/3}Z_{\mathrm{es}}^N(r)\,dr\right)\\
&\quad - N^{-2/3}Y_1\left(\int_0^t N^{5/3}Z_{\mathrm{e}}^N(r)Z_{\mathrm{s}}^N(r)\,dr\right),\\
Z_{\mathrm{es}}^N(t) &= -N^{-2/3}Y_{-1}\left(\int_0^t N^{2/3}Z_{\mathrm{es}}^N(r)\,dr\right)\\
&\quad + N^{-2/3}Y_1\left(\int_0^t N^{5/3}Z_{\mathrm{e}}^N(r)Z_{\mathrm{s}}^N(r)\,dr\right)\\
&\quad - N^{-2/3}Y_2\left(\int_0^t N^{1/3}Z_{\mathrm{es}}^N(r)\,dr\right),\\
Z_{\mathrm{e}}^N(t) &= 1 + N^{-1}Y_{-1}\left(\int_0^t N^{2/3}Z_{\mathrm{es}}^N(r)\,dr\right)\\
&\quad - N^{-1}Y_1\left(\int_0^t N^{5/3}Z_{\mathrm{e}}^N(r)Z_{\mathrm{s}}^N(r)\,dr\right)\\
&\quad + N^{-1}Y_2\left(\int_0^t N^{1/3}Z_{\mathrm{es}}^N(r)\,dr\right),\\
Z_{\mathrm{p}}^N(t) &= N^{-2/3}Y_2\left(\int_0^t N^{1/3}Z_{\mathrm{es}}^N(r)\,dr\right).
\end{aligned}$$

Rescaling time by $N^{1/3}$ and defining $V_i^N(t) = Z_i^N(N^{1/3}t)$, we have

$$\begin{aligned}
V_{\mathrm{s}}^N(t) &= 1 + N^{-2/3}Y_{-1}\left(\int_0^t NV_{\mathrm{es}}^N(r)\,dr\right)\\
&\quad - N^{-2/3}Y_1\left(\int_0^t N^2V_{\mathrm{e}}^N(r)V_{\mathrm{s}}^N(r)\,dr\right),\\
V_{\mathrm{es}}^N(t) &= -N^{-2/3}Y_{-1}\left(\int_0^t NV_{\mathrm{es}}^N(r)\,dr\right)
\end{aligned}$$

$$
\begin{aligned}
&+ N^{-2/3}Y_1\left(\int_0^t N^2 V_{\mathrm{e}}^N(r)V_{\mathrm{s}}^N(r)\,dr\right)\\
&- N^{-2/3}Y_2\left(\int_0^t N^{2/3} V_{\mathrm{es}}^N(r)\,dr\right),\\
V_{\mathrm{e}}^N(t) &= 1 + N^{-1}Y_{-1}\left(\int_0^t N V_{\mathrm{es}}^N(r)\,dr\right)\\
&- N^{-1}Y_1\left(\int_0^t N^2 V_{\mathrm{e}}^N(r)V_{\mathrm{s}}^N(r)\,dr\right)\\
&+ N^{-1}Y_2\left(\int_0^t N^{2/3} V_{\mathrm{es}}^N(r)\,dr\right),\\
V_{\mathrm{p}}^N(t) &= N^{-2/3}Y_2\left(\int_0^t N^{2/3} V_{\mathrm{es}}^N(r)\,dr\right),\\
V_{\mathrm{s}}^N(t) + V_{\mathrm{es}}^N(t) &= 1 - N^{-2/3}Y_2\left(\int_0^t N^{2/3} V_{\mathrm{es}}^N(r)\,dr\right).
\end{aligned}
$$

Note that $V_{\mathrm{s}}^N + V_{\mathrm{es}}^N \le 1$ and $V_{\mathrm{e}}^N + N^{-1/3}V_{\mathrm{es}}^N = 1$, so $\sup_{r\le t}|V_{\mathrm{e}}^N(r) - 1| \to 0$. It follows that for $0 < \rho < 1$ and $N$ sufficiently large, $\inf_{r\le t} V_{\mathrm{e}}^N(r) \ge \rho$ and $V_{\mathrm{s}}^N(t) \le \widehat{V}_{\mathrm{s}}^N(t)$, where $\widehat{V}_{\mathrm{s}}^N$ is the solution of

$$
\widehat{V}_{\mathrm{s}}^N(t) = 1 + N^{-2/3}Y_{-1}\left(\int_0^t N V_{\mathrm{es}}^N(r)\,dr\right) - N^{-2/3}Y_1\left(\int_0^t N^2\rho \widehat{V}_{\mathrm{s}}^N(r)\,dr\right).
$$

The fact that $\sup_{r\le t}\frac{N V_{\mathrm{es}}^N(r)}{N^2\rho} \to 0$ ensures that $\sup_{\delta<r\le t} V_{\mathrm{s}}^N(r) \to 0$, for $0 < \delta < t$, and $(V_{\mathrm{p}}^N, V_{\mathrm{s}}^N + V_{\mathrm{es}}^N)$ converges to the solution of

$$
\begin{aligned}
V_{\mathrm{p}}(t) &= \int_0^t V_{\mathrm{es}}(r)\,dr,\\
V_{\mathrm{es}}(t) &= 1 - \int_0^t V_{\mathrm{es}}(r)\,dr,
\end{aligned}
$$

that is, $V_{\mathrm{es}}(t) = e^{-t}$ and $V_{\mathrm{p}}(t) = 1 - e^{-t}$. On the original time scale $Z_{\mathrm{p}}^N(t) \approx 1 - e^{-t/10}$, that is, $X_{\mathrm{p}}(t) \approx 100(1 - e^{-t/10})$, which matches well the simulation results in the lower plot in Figure 1 of [9].

Rao and Arkin also consider $X_{\mathrm{s}}(0) = 100$, $X_{\mathrm{e}}(0) = 10$, $\kappa_1 = \kappa_{-1} = 1$, $\kappa_2 = 0.1$. For this example, let $N = 100$ and define $Z_{\mathrm{s}}^N = N^{-1}X_{\mathrm{s}}$, $Z_{\mathrm{es}}^N = N^{-1/2}X_{\mathrm{es}}$, $Z_{\mathrm{e}}^N = N^{-1/2}X_{\mathrm{e}}$, and $Z_{\mathrm{p}}^N = X_{\mathrm{p}}$, and set $\kappa_2 = N^{-1/2}$. Then the normalized system becomes

$$
\begin{aligned}
Z_{\mathrm{s}}^N(t) &= 1 + N^{-1}Y_{-1}\left(\int_0^t N^{1/2} Z_{\mathrm{es}}^N(r)\,dr\right) - N^{-1}Y_1\left(\int_0^t N^{3/2} Z_{\mathrm{e}}^N(r)Z_{\mathrm{s}}^N(r)\,dr\right),\\
Z_{\mathrm{es}}^N(t) &= -N^{-1/2}Y_{-1}\left(\int_0^t N^{1/2} Z_{\mathrm{es}}^N(r)\,dr\right) + N^{-1/2}Y_1\left(\int_0^t N^{3/2} Z_{\mathrm{e}}^N(r)Z_{\mathrm{s}}^N(r)\,dr\right)
\end{aligned}
$$

$$- N^{-1/2}Y_2\Big(\int_0^t Z_{\rm es}^N(r)\,dr\Big),$$

$$Z_{\rm e}^N(t) = 1 + N^{-1/2}Y_{-1}\Big(\int_0^t N^{1/2}Z_{\rm es}^N(r)\,dr\Big)$$
$$- N^{-1/2}Y_1\Big(\int_0^t N^{3/2}Z_{\rm e}^N(r)Z_{\rm s}^N(r)\,dr\Big) + N^{-1/2}Y_2\Big(\int_0^t Z_{\rm es}^N(r)\,dr\Big),$$

$$Z_{\rm p}^N(t) = Y_2\Big(\int_0^t Z_{\rm es}^N(r)\,dr\Big).$$

Since

$$N^{-1/2}Y_1\Big(\int_0^t N^{3/2}Z_{\rm e}^N(r)Z_{\rm s}^N(r)\,dr\Big)$$
$$\le 1 + N^{-1/2}Y_{-1}\Big(\int_0^t N^{1/2}Z_{\rm es}^N(r)\,dr\Big)$$
$$+ N^{-1/2}Y_2\Big(\int_0^t Z_{\rm es}^N(r)\,dr\Big)$$

and $Z_{\rm es}(t) \le 1$, it follows that $\sup_{r\le t}|Z_{\rm s}^N(r) - 1| \to 0$, for each $t > 0$. Noting that $\int_0^t Z_{\rm es}^N(r)\,dr = t - \int_0^t Z_{\rm e}^N(r)\,dr$, we must have $\int_0^t Z_{\rm es}^N(r)\,dr \approx t$ and $Z_{\rm p}^N(t) \approx Y_2(t)$.

2.3. *Reversible isomerization.* Next, we consider a model of reversible isomerization studied by Cao, Gillespie and Petzold [3]. The model involves three chemical species and two reactions:

$$S_1 \underset{\kappa_2}{\overset{\kappa_1}{\rightleftarrows}} S_2, \qquad S_2 \xrightarrow{\kappa_3} S_3.$$

The model is given by

$$X_1(t) = X_1(0) - Y_1\Big(\int_0^t \kappa_1 X_1(s)\,ds\Big) + Y_2\Big(\int_0^t \kappa_2 X_2(s)\,ds\Big),$$
$$X_2(t) = X_2(0) + Y_1\Big(\int_0^t \kappa_1 X_1(s)\,ds\Big)$$
$$- Y_2\Big(\int_0^t \kappa_2 X_2(s)\,ds\Big) - Y_3\Big(\int_0^t \kappa_3 X_2(s)\,ds\Big),$$
$$X_3(t) = X_3(0) + Y_3\Big(\int_0^t \kappa_3 X_2(s)\,ds\Big).$$

The first set of parameter values in (34) of [3] give $X_1(0) = 1200$, $X_2(0) = 600$, $X_3(0) = 0$ and $\kappa_1 = 1$, $\kappa_2 = 2$, $\kappa_3 = 5 \times 10^{-5}$. Let $N = 1000$, and define $Z_1^N = N^{-1}X_1, Z_2^N = N^{-1}X_2, Z_3^N = X_3$, and $\kappa_3 = 5N^{-5/3}$. Then the normalized system becomes

$$\begin{aligned}
Z_1^N(t) &= Z_1^N(0) - N^{-1}Y_1\biggl(\int_0^t NZ_1^N(s)\,ds\biggr) + N^{-1}Y_2\biggl(\int_0^t 2NZ_2^N(s)\,ds\biggr),\\
Z_2^N(t) &= Z_2^N(0) + N^{-1}Y_1\biggl(\int_0^t NZ_1^N(s)\,ds\biggr)\\
&\qquad - N^{-1}Y_2\biggl(\int_0^t 2NZ_2^N(s)\,ds\biggr) - N^{-1}Y_3\biggl(\int_0^t 5N^{-2/3}Z_2^N(s)\,ds\biggr),\\
Z_3^N(t) &= Y_3\biggl(\int_0^t 5N^{-2/3}Z_2^N(s)\,ds\biggr).
\end{aligned}$$

Assuming that $(Z_1^N(0), Z_2^N(0)) \to (Z_1(0), Z_2(0))$ (which gives $Z_1(0) = 1.2$ and $Z_2(0) = 0.6$ for the particular values in [3]), the limiting system is

$$\begin{aligned}
Z_1(t) &= Z_1(0) - \int_0^t Z_1(s)\,ds + \int_0^t 2Z_2(s)\,ds,\\
Z_2(t) &= Z_2(0) + \int_0^t Z_1(s)\,ds - \int_0^t 2Z_2(s)\,ds,\\
Z_3(t) &= 0.
\end{aligned}$$

Consequently,

$$Z_1(t) + Z_2(t) = Z_1(0) + Z_2(0),$$

$$D(t) \equiv Z_1(t) - 2Z_2(t) = D(0) - 3\int_0^t D(s)\,ds,$$

so

$$\begin{aligned}
Z_1(t) &= \tfrac{1}{3}(Z_1(0) - 2Z_2(0))e^{-3t} + \tfrac{2}{3}(Z_1(0) + Z_2(0)),\\
Z_2(t) &= -\tfrac{1}{3}(Z_1(0) - 2Z_2(0))e^{-3t} + \tfrac{1}{3}(Z_1(0) + Z_2(0))
\end{aligned}$$

and

$$\lim_{t\to\infty}(Z_1(t), Z_2(t)) = (\tfrac{2}{3}, \tfrac{1}{3})(Z_1(0) + Z_2(0)).$$

Defining $U_i^N(t) = Z_i^N(N^{2/3}t)$, the system becomes

$$\begin{aligned}
U_1^N(t) &= U_1^N(0) - N^{-1}Y_1\biggl(\int_0^t N^{5/3}U_1^N(s)\,ds\biggr) + N^{-1}Y_2\biggl(\int_0^t 2N^{5/3}U_2^N(s)\,ds\biggr),\\
U_2^N(t) &= U_2^N(0) + N^{-1}Y_1\biggl(\int_0^t N^{5/3}U_1^N(s)\,ds\biggr) - N^{-1}Y_2\biggl(\int_0^t 2N^{5/3}U_2^N(s)\,ds\biggr)\\
&\qquad - N^{-1}Y_3\biggl(\int_0^t 5U_2^N(s)\,ds\biggr),\\
U_3^N(t) &= Y_3\biggl(\int_0^t 5U_2^N(s)\,ds\biggr)
\end{aligned}$$

and, hence,

$$U_1^N(t) + U_2^N(t) = U_1^N(0) + U_2^N(0) - N^{-1}Y_3\biggl(\int_0^t 5U_2^N(s)\,ds\biggr).$$

Dividing the equation for $U_1^N$ by $N^{2/3}$, it follows that

$$\lim_{N\to\infty}\biggl(\int_0^t U_1^N(s)\,ds - \int_0^t 2U_2^N(s)\,ds\biggr) = 0$$

and, hence, assuming $\lim_{N\to\infty}(U_1^N(0) + U_2^N(0)) = C$,

$$\lim_{N\to\infty}\int_0^t U_2^N(s)\,ds = \tfrac{1}{3}Ct$$

and $U_3^N$ converges to

$$U_3(t) = Y_3\biggl(\frac{5Ct}{3}\biggr),$$

that is, a Poisson process with parameter $5C/3$.

Rescaling time by $N^{5/3}$, and defining $V_i^N(t) = N^{-1}X_i(N^{5/3}t)$, we have

$$\begin{aligned}
V_1^N(t) &= V_1^N(0) - N^{-1}Y_1\biggl(\int_0^t N^{8/3}V_1^N(s)\,ds\biggr) + N^{-1}Y_2\biggl(\int_0^t 2N^{8/3}V_2^N(s)\,ds\biggr),\\
V_2^N(t) &= V_2^N(0) + N^{-1}Y_1\biggl(\int_0^t N^{8/3}V_1^N(s)\,ds\biggr) - N^{-1}Y_2\biggl(\int_0^t 2N^{8/3}V_2^N(s)\,ds\biggr)\\
&\qquad - N^{-1}Y_3\biggl(\int_0^t 5NV_2^N(s)\,ds\biggr),\\
V_3^N(t) &= N^{-1}Y_3\biggl(\int_0^t 5NV_2^N(s)\,ds\biggr).
\end{aligned}$$

Setting $R^N(t) = V_1^N(t) + V_2^N(t)$ and assuming $\lim_{N\to\infty} R^N(0) = R(0)$, we have $(R^N, V_3^N) \to (R, V_3)$ satisfying

$$R(t) = R(0) - \int_0^t \tfrac{5}{3}R(s)\,ds, \qquad V_3(t) = \int_0^t \tfrac{5}{3}R(s)\,ds,$$

which gives

$$R(t) = R(0)e^{-(5/3)t}, \qquad V_3(t) = R(0)(1 - e^{-(5/3)t}).$$

For $\delta > 0$,

$$\lim_{N\to\infty}\sup_{\delta\le r\le t}(|V_1^N(r) - \tfrac{2}{3}R(t)| + |V_2^N(r) - \tfrac{1}{3}R(t)|) = 0.$$

Note that the simulation results given in Figure 1 of [3] appear to be plots of $\{(X_3(t_k)\delta, X_i(t_k))\}$, for some $\delta > 0$, rather than of $\{(t_k, X_i(t_k))\}$, where the $\{t_k\}$ are the jump times of $X_3$. Consequently, their plots show linear decay rather than exponential decay.

Cao, Gillespie and Petzold also study a second set of parameter values (35) in [3], taking $X_1(0) = 2000, X_2(0) = X_3(0) = 0$ and $\kappa_1 = 10, \kappa_2 = 4 \times 10^4$, $\kappa_3 = 2$. Letting $N = 10^4$, we now define $Z_1^N = N^{-1}X_1, Z_2^N = X_2, Z_3^N = X_3$. The normalized system becomes

$$Z_1^N(t) = 0.2 - N^{-1}Y_1\left(\int_0^t 10NZ_1^N(s)\,ds\right) + N^{-1}Y_2\left(\int_0^t 4NZ_2^N(s)\,ds\right),$$

$$Z_2^N(t) = Y_1\left(\int_0^t 10NZ_1^N(s)\,ds\right) - Y_2\left(\int_0^t 4NZ_2^N(s)\,ds\right) - Y_3\left(\int_0^t 2Z_2^N(s)\,ds\right),$$

$$Z_3^N(t) = Y_3\left(\int_0^t 2Z_2^N(s)\,ds\right).$$

Dividing the equation for $Z_2^N$ by $N$, we see that

$$\sup_{r\le t}\frac{1}{N}\left(Y_1\left(\int_0^r 10NZ_1^N(s)\,ds\right) - Y_2\left(\int_0^r 4NZ_2^N(s)\,ds\right)\right) \to 0,$$

for each $t > 0$, and hence,

$$\sup_{r\le t}|Z_1^N(r) - 0.2| \to 0,$$

$$\sup_{r\le t}\left|\int_0^r Z_2^N(s)\,ds - \tfrac{1}{2}r\right| \to 0,$$

and $Z_3^N$ converges to $Y_3(t)$.

Rescaling time by $N$, and defining $V_i^N(t) = N^{-1}X_i(Nt)$, for $i = 1,3$ and $V_2^N(t) = X_2(Nt)$, we have

$$V_1^N(t) = 0.2 - N^{-1}Y_1\left(\int_0^t 10N^2V_1^N(s)\,ds\right) + N^{-1}Y_2\left(\int_0^t 4N^2V_2^N(s)\,ds\right),$$

$$V_2^N(t) = Y_1\left(\int_0^t 10N^2V_1^N(s)\,ds\right) - Y_2\left(\int_0^t 4N^2V_2^N(s)\,ds\right) - Y_3\left(\int_0^t 2NV_2^N(s)\,ds\right),$$

$$V_3^N(t) = N^{-1}Y_3\left(\int_0^t 2NV_2^N(s)\,ds\right).$$

Let $\widehat{V}_1^N(t) = V_1^N(t) + N^{-1}V_2^N(t)$, and since earlier results imply $V_2^N(t) = Z_2^N(Nt) \approx \frac{1}{2}$, we have $\sup_{r\le t}|V_1^N(r) - \widehat{V}_1^N(r)| \to 0$.
Then

$$\widehat{V}_1^N(t) = 0.2 - N^{-1}Y_3\left(\int_0^t 2NV_2^N(s)\,ds\right),$$

and from the equation for $V_2^N$,

$$\sup_{r\le t}\left|\int_0^r 10\widehat{V}_1^N(s)\,ds - \int_0^r (4+2N^{-1})V_2^N(s)\right| \to 0.$$

Consequently, $(V_1^N, V_3^N)$ converge to the solution of

$$V_1(t) = 0.2 - \int_0^t 5V_1(s)\,ds,$$
$$V_3(t) = \int_0^t 5V_1(s)\,ds,$$

giving

$$V_1(t) = 0.2e^{-5t}, \qquad V_3(t) = 0.2(1-e^{-5t}).$$

To better understand the behavior of $V_2^N$, for a bounded function $f(v_2)$, define

$$A^N f(v_1,v_2) = 10v_1(f(v_2+1)-f(v_2)) + (4+2N^{-1})v_2(f(v_2-1)-f(v_2))$$

and note that

$$f(V_2^N(t)) - N^2\int_0^t A^N f(V_1^N(s), V_2^N(s))\,ds$$

is a martingale. Dividing by $N^2$, it follows that

$$\int_0^t A^N f(V_1^N(s), V_2^N(s))\,ds \to 0$$

and that for each $t$, $V_2^N(t)$ converges in distribution to a random variable $V_2(t)$ satisfying

$$E[Af(V_1(t), V_2(t))] = 0,$$

where

$$Af(v_1,v_2) = 10v_1(f(v_2+1)-f(v_2)) + 4v_2(f(v_2-1)-f(v_2))$$

(see, e.g., [8]).

For each $v_1$, $A_{v_1}f(v_2) \equiv Af(v_1,v_2)$ is the generator of an infinite-server queueing model with arrival rate $10v_1$ and service rate 4. It follows that $V_2(t)$ has a Poisson distribution with parameter $2.5V_1(t)$. Note that $V_2^N$ does not converge in a functional sense. In particular, for $0 < t_1 < t_2 < \cdots < t_m$, $(V_2^N(t_1), V_2^N(t_2), \ldots, V_2^N(t_m))$ converges in distribution and the components of the limit $(V_2(t_1), V_2(t_2), \ldots, V_2(t_m))$ are independent Poisson random variables.

**3. Intracellular viral kinetics.** Next we consider a model of an intracellular viral infection given in [10] and studied further in [6]. We follow the presentation (in particular, the indexing) in [6]. The model includes three time-varying species, the viral template, the viral genome and the viral structural protein. We denote these as species 1, 2 and 3, respectively, and let $X_i(t)$ denote the number of molecules of species $i$ in the system at time $t$. The model involves six reactions, designated (28a)–(28f) in [6]:

(a) $T + \text{stuff} \xrightarrow{\kappa_1} T + G$,
(b) $G \xrightarrow{\kappa_2} T$,
(c) $T + \text{stuff} \xrightarrow{\kappa_3} T + S$,
(d) $T \xrightarrow{\kappa_4} \varnothing$,
(e) $S \xrightarrow{\kappa_5} \varnothing$,
(f) $G + S \xrightarrow{\kappa_6} V$,

where "stuff" refers to nucleotides and amino acids that are assumed available at constant concentrations. The basic model satisfies

$$\begin{aligned} X_1(t) &= X_1(0) + Y_b\left(\int_0^t \kappa_2 X_2(s)\,ds\right) - Y_d\left(\int_0^t \kappa_4 X_1(s)\,ds\right)\\ X_2(t) &= X_2(0) + Y_a\left(\int_0^t \kappa_1 X_1(s)\,ds\right) - Y_b\left(\int_0^t \kappa_2 X_2(s)\,ds\right)\\ &\quad - Y_f\left(\int_0^t \kappa_6 X_2(s) X_3(s)\,ds\right)\\ X_3(t) &= X_3(0) + Y_c\left(\int_0^t \kappa_3 X_1(s)\,ds\right) - Y_e\left(\int_0^t \kappa_5 X_3(s)\,ds\right)\\ &\quad - Y_f\left(\int_0^t \kappa_6 X_2(s) X_3(s)\,ds\right). \end{aligned}$$

Following Haseltine and Rawlings, $X_1(0) = 1, X_2(0) = X_3(0) = 0$, while the reaction constants from Table III of [6] are as given below. Let $N = 1000$, corresponding to the order of magnitude of the largest reaction constant. Then the rate constants can be expressed as follows:

| | | |
|---|---|---|
| $\kappa_1$ | 1 | 1 |
| $\kappa_2$ | 0.025 | $2.5N^{-2/3}$ |
| $\kappa_3$ | 1000 | $N$ |
| $\kappa_4$ | 0.25 | 0.25 |
| $\kappa_5$ | 2 | 2 |
| $\kappa_6$ | $7.5 \times 10^{-6}$ | $0.75N^{-5/3}$. |

Note that, for simplicity, we have replaced $\kappa_5 = 1.9985$ by $\kappa_5 = 2$.

We have identified $N$ with the largest rate constant, but it is also the order of magnitude of the most abundent species. Writing $Z_1^N = X_1$, $Z_2^N =$

$N^{-2/3}X_2$ and $Z_3^N = N^{-1}X_3$, the normalized system with the scaled rate constants becomes

$$\begin{aligned}
Z_1^N(t) &= Z_1^N(0) + Y_b\Big(\int_0^t 2.5 Z_2^N(s)\,ds\Big) - Y_d\Big(\int_0^t 0.25 Z_1^N(s)\,ds\Big),\\
Z_2^N(t) &= Z_2^N(0) + N^{-2/3}Y_a\Big(\int_0^t Z_1^N(s)\,ds\Big) - N^{-2/3}Y_b\Big(\int_0^t 2.5 Z_2^N(s)\,ds\Big)\\
&\qquad - N^{-2/3}Y_f\Big(\int_0^t 0.75 Z_2^N(s)Z_3^N(s)\,ds\Big),\\
Z_3^N(t) &= Z_3^N(0) + N^{-1}Y_c\Big(\int_0^t N Z_1^N(s)\,ds\Big) - N^{-1}Y_e\Big(\int_0^t 2N Z_3^N(s)\,ds\Big)\\
&\qquad - N^{-1}Y_f\Big(\int_0^t 0.75 Z_2^N(s)Z_3^N(s)\,ds\Big).
\end{aligned}$$

We also write $X_i^N$ for the system $X_i$ with rate constants expressed in terms of $N$, including the superscript $N$ only to emphasize the dependence of the model on the scaling parameter.

There is substantial probability that the infection dies out quickly, but if $Z_2^N$ reaches any significant level, the chance becomes negligible.

To be precise, let

$$K_a^N(t) = Y_a\Big(\int_0^t X_1^N(s)\,ds\Big),$$

and for $k = 1, 2, \ldots,$ define

$$\beta_k^N = \inf\{t \geq 0 : K_a^N(t) \geq k\}. \tag{3.1}$$

We have the following results.

Lemma 3.1 (Probability of infection). *For $\beta_k^N$ defined in (3.1),*

$$\lim_{k\to\infty}\lim_{N\to\infty} P\{\beta_k^N < \infty\} = 0.75.$$

*In particular, there exist $k_N \to \infty$ such that $\lim_{N\to\infty} P\{\beta_{k_N}^N < \infty\} = 0.75$.*

Remark 3.2. The theorem essentially gives the probability that a single virus successfully infects the cell. The argument, which essentially compares the initial stages of the infection to a branching process, is a standard tool in the analysis of epidemic models. See, for example, [2].

Proof of Lemma 3.1. To understand the initial behavior of the system, consider

$$X_1^N(t) = Y_b\Big(\int_0^t 2.5N^{-2/3}X_2^N(s)\,ds\Big) - Y_d\Big(\int_0^t 0.25 X_1^N(s)\,ds\Big), \tag{3.2}$$

$$
\begin{aligned}
X_2^N(t) &= 1 + Y_a\left(\int_0^t X_1^N(s)\,ds\right) - Y_b\left(\int_0^t 2.5N^{-2/3}X_2^N(s)\,ds\right) \\
&\quad - Y_f\left(\int_0^t 0.75N^{-2/3}X_2^N(s)Z_3^N(s)\,ds\right),
\end{aligned} \tag{3.3}
$$

$$
\begin{aligned}
Z_3^N(t) &= N^{-1}Y_c\left(\int_0^t NX_1^N(s)\,ds\right) - N^{-1}Y_e\left(\int_0^t 2NZ_3^N(s)\,ds\right) \\
&\quad - N^{-1}Y_f\left(\int_0^t 0.75N^{-2/3}X_2^N(s)Z_3^N(s)\,ds\right).
\end{aligned} \tag{3.4}
$$

If the virus production term (reaction $f$) is dropped from the equation for $X_2^N$, the resulting system,

$$
\begin{aligned}
\widehat{X}_1^N(t) &= Y_b\left(\int_0^t 2.5N^{-2/3}\widehat{X}_2^N(s)\,ds\right) - Y_d\left(\int_0^t 0.25\widehat{X}_1^N(s)\,ds\right), \\
\widehat{X}_2^N(t) &= 1 + Y_a\left(\int_0^t \widehat{X}_1^N(s)\,ds\right) - Y_b\left(\int_0^t 2.5N^{-2/3}\widehat{X}_2^N(s)\,ds\right),
\end{aligned}
$$

determines a continuous time, two-type branching process. It is easy to check that the process is supercritical. The "lifetime" of each Type 1 molecule is exponentially distributed with parameter 0.25 and the number $\xi_i$ of Type 2 molecules created by the $i$th Type 1 molecule during its lifetime has a shifted geometric distribution with expectation 4, that is,

$$
P\{\xi = k\} = \int_0^\infty 0.25e^{-0.25t}e^{-t}\frac{t^k}{k!}\,dt = \frac{1}{5}\left(\frac{4}{5}\right)^k, \qquad k = 0, 1, \ldots.
$$

Each Type 2 molecule is eventually converted to Type 1. Starting with a single Type 1 or Type 2 molecule, the probability of extinction is simply the probability that

$$
S_n = 1 + \sum_{i=1}^n (\xi_i - 1)
$$

hits zero for some $n \ge 0$, an event with probability 0.25.

To complete the proof of the lemma, one only needs to check that

$$
Y_f\left(\int_0^{\beta_k^N} 0.75N^{-2/3}X_2^N(s)Z_3^N(s)\,ds\right) \to 0
$$

in probability for each $k$. But

$$
\begin{aligned}
&E\left[Y_f\left(\int_0^{\beta_k^N} 0.75N^{-2/3}X_2^N(s)Z_3^N(s)\,ds\right)\right] \\
&\qquad \le 0.75N^{-2/3}(k+1)E\left[\int_0^{\beta_k^N} Z_3^N(s)\,ds\right]
\end{aligned}
$$

$$\le 0.375N^{-2/3}(k+1)E\left[\int_0^{\beta_k^N} X_1^N(s)\,ds\right]$$
$$\le 0.375N^{-2/3}(k+1)k$$
$$\to 0,$$

where the second inequality follows from equation (3.4) and the last inequality follows from the fact that

$$E\left[\int_0^{\beta_k^N} X_1^N(s)\,ds\right] = E\left[Y_a\left(\int_0^{\beta_k^N} X_1^N(s)\,ds\right)\right] \le k. \qquad \square$$

We now want to describe the behavior of the process once the infection is established. Since we have scaled $X_2^N$ by $N^{-2/3}$, $X_2^N$ must reach a level that is $O(N^{2/3})$ to have nontrivial behavior.

As in the proof of Lemma 3.1, if the virus production term is dropped from the equation for $X_2^N$, the expectation of the resulting two-type branching process $(\widehat{X}_1^N, \widehat{X}_2^N)$ satisfies $\dot{m}(t) = Q^N m(t)$, where

$$Q^N = \begin{pmatrix} -0.25 & 2.5N^{-2/3} \\ 1 & -2.5N^{-2/3} \end{pmatrix}.$$

Following the classical analysis for branching processes (see Section V.7 of [1]), the largest eigenvalue for $Q^N$ is

$$\lambda^N = \frac{-(0.25 + 2.5N^{-2/3}) + \sqrt{(0.25 + 2.5N^{-2/3})^2 + 7.5N^{-2/3}}}{2}$$

and the total "population" should grow like $e^{\lambda^N t}$. There exists $\rho^N > 0$ satisfying

$$(1 - 0.25\rho^N) = \lambda^N \rho^N, \qquad 2.5(\rho^N - 1)N^{-2/3} = \lambda^N,$$

that is, $(\rho^N, 1)$ is the corresponding left eigenvector. It follows that $\rho^N \to 4$ and $N^{2/3}\lambda^N \to 7.5$.

Let $R^N(t) = \rho^N X_1^N(t) + X_2^N(t)$, and define

$$\tau_\varepsilon^N = \inf\{t : R^N(t) \ge N^{2/3}\varepsilon\}.$$

We are really interested in the first time $X_2^N$ reaches $N^{2/3}\varepsilon$, but defining $\tau_\varepsilon^N$ in terms of $R^N$ rather than $X_2^N$ simplifies the proof of the next theorem.

THEOREM 3.3 (Time until establishment). *For $k_N$ as in Lemma 3.1 and each $0 < \varepsilon < 2$ and $\delta > 0$,*

$$\lim_{N\to\infty} P\left\{\left|\frac{\tau_\varepsilon^N}{N^{2/3}\log N} - \frac{4}{45}\right| > \delta \,\middle|\, \beta_{k_N}^N < \infty\right\} = 0.$$

*In particular,* $\lim_{N\to\infty} P\{\tau_\varepsilon^N < \infty\} = 0.75$.

PROOF. In the calculations that follow, recall that the law of large numbers implies that, for a unit Poisson process,

$$\lim_{u_0\to\infty}\sup_{u\ge u_0}\left|\frac{Y(u)}{u}-1\right|=0 \qquad \text{a.s.} \tag{3.5}$$

In addition, note that, without loss of generality, we can assume that $k_N/\log N\to 0$.

Define

$$\begin{aligned}
K_a^N(t)&=Y_a\biggl(\int_0^t X_1^N(s)\,ds\biggr),\\
\tilde K_a^N(t)&=Y_a\biggl(\int_0^t X_1^N(s)\,ds\biggr)-\int_0^t X_1^N(s)\,ds,\\
K_b^N(t)&=Y_b\biggl(\int_0^t 2.5N^{-2/3}X_2^N(s)\,ds\biggr),\\
\tilde K_b^N(t)&=Y_b\biggl(\int_0^t 2.5N^{-2/3}X_2^N(s)\,ds\biggr)-\int_0^t 2.5N^{-2/3}X_2^N(s)\,ds,
\end{aligned}$$

and similarly for $K_c^N$, $K_d^N$ and so on. Since $\int_0^{\beta_{k_N}^N} X_1^N(s)\,ds$ is the $k_N$th jump time of $Y_a$,

$$\mathbf{1}_{\{\beta_{k_N}^N<\infty\}}\left|\frac{1}{k_N}\int_0^{\beta_{k_N}^N} X_1^N(s)\,ds-1\right|\to 0,$$

and it follows from (3.5) that

$$\mathbf{1}_{\{\beta_{k_N}^N<\infty\}}\sup_{t\ge\beta_{k_N}^N}\left|\frac{K_a^N(t)}{\int_0^t X_1^N(s)\,ds}-1\right|\to 0,$$

$$\mathbf{1}_{\{\beta_{k_N}^N<\infty\}}\sup_{t\ge\beta_{k_N}^N}\left|\frac{K_b^N(t)}{\int_0^t 2.5N^{-2/3}X_2^N(s)\,ds}-1\right|\to 0$$

and

$$\mathbf{1}_{\{\beta_{k_N}^N<\infty\}}\sup_{t\ge\beta_{k_N}^N}\left|\frac{K_d^N(t)}{K_a^N(t)}-0.25\right|\to 0.$$

With reference to (3.4), $K_e^N(t)\le K_c^N(t)$ and

$$\begin{aligned}
&\lim_{N\to\infty}\mathbf{1}_{\{\beta_{k_N}^N<\infty\}}\sup_{t\ge\beta_{k_N}^N}\biggl(\frac{\int_0^t 2Z_3^N(s)\,ds}{\int_0^t X_1^N(s)\,ds}-1\biggr)\vee 0\\
&\qquad=\lim_{N\to\infty}\mathbf{1}_{\{\beta_{k_N}^N<\infty\}}\sup_{t\ge\beta_{k_N}^N}\biggl(\frac{K_e^N(t)}{K_c^N(t)}-1\biggr)\vee 0=0.
\end{aligned}$$

Consequently,

$$\begin{aligned}
&\lim_{N\to\infty} \mathbf{1}_{\{\beta^N_{k_N}<\infty\}} \sup_{\beta^N_{k_N}\le t\le \tau^N_\varepsilon} \left(\frac{K^N_f(t)}{1+K^N_a(t)} - 0.375\varepsilon\right)\vee 0\\
&\quad\le \lim_{N\to\infty} \mathbf{1}_{\{\beta^N_{k_N}<\infty\}} \sup_{\beta^N_{k_N}\le t\le \tau^N_\varepsilon} \left(\frac{Y_f(\int_0^t 0.75\varepsilon Z^N_3(s)\,ds)}{1+K^N_a(t)} - 0.375\varepsilon\right)\vee 0\\
&\quad\le \lim_{N\to\infty} \mathbf{1}_{\{\beta^N_{k_N}<\infty\}} \sup_{\beta^N_{k_N}\le t\le \tau^N_\varepsilon} \left(\frac{\int_0^t 0.75\varepsilon Z^N_3(s)\,ds}{\int_0^t X^N_1(s)\,ds} - 0.375\varepsilon\right)\vee 0\\
&\quad= 0.
\end{aligned}$$

Since $K^N_d(t)\le K^N_b(t)\le 1+K^N_a(t)$ and

$$\begin{aligned}
R^N(t) &= 1+K^N_a(t) - \rho^N K^N_d(t) + (\rho^N-1)K^N_b(t) - K^N_f(t)\\
&\ge 1+K^N_a(t) - K^N_d(t) - K^N_f(t),
\end{aligned}$$

we have

$$\begin{aligned}
\frac{R^N(t)}{K^N_d(t)} &\ge \frac{R^N(t)}{K^N_b(t)}\\
&\ge \frac{R^N(t)}{1+K^N_a(t)} \qquad (3.6)\\
&\ge 1-\frac{K^N_d(t)+K^N_f(t)}{1+K^N_a(t)}
\end{aligned}$$

and

$$\begin{aligned}
&\lim_{N\to\infty} \mathbf{1}_{\{\beta^N_{k_N}<\infty\}} \inf_{\beta^N_{k_N}\le t\le \tau^N_\varepsilon} \left(1-\frac{K^N_d(t)+K^N_f(t)}{1+K^N_a(t)} - 0.75(1-0.5\varepsilon)\right)\wedge 0\\
&\quad= 0. \qquad (3.7)
\end{aligned}$$

In other words, on the event $\{\beta^N_{k_N}<\infty\}$,

$$\inf_{\beta^N_{k_N}\le t\le \tau^N_\varepsilon}\left(1-\frac{K^N_d(t)+K^N_f(t)}{1+K^N_a(t)}\right)$$

is asymptotically bounded below by $0.75(1-0.5\varepsilon)>0$.

Writing $\rho$ instead of $\rho^N$, and $\beta$ instead of $\beta^N_{k_N}$,

$$\log R^N(t) = \log R^N(\beta^N_{k_N}) + \int_\beta^t \log\frac{R^N(s-)+1}{R^N(s-)}\,dK^N_a(s)$$

$$
\begin{aligned}
&+\int_\beta^t \log\frac{R^N(s-)+(\rho-1)}{R^N(s-)}\,dK_b^N(s)\\
&+\int_\beta^t \log\frac{R^N(s-)-\rho}{R^N(s-)}\,dK_d^N(s)+\int_\beta^t \log\frac{R^N(s-)-1}{R^N(s-)}\,dK_f^N(s)\\
=&\log R^N(\beta_{k_N}^N)+\int_\beta^t\left(\log\frac{R^N(s-)+1}{R^N(s-)}-\frac{1}{R^N(s-)}\right)dK_a^N(s)\\
&+\int_\beta^t\left(\log\frac{R^N(s-)+(\rho-1)}{R^N(s-)}-\frac{\rho-1}{R^N(s-)}\right)dK_b^N(s)\\
&+\int_\beta^t\left(\log\frac{R^N(s-)-\rho}{R^N(s-)}+\frac{\rho}{R^N(s-)}\right)dK_d^N(s)\\
&+\int_\beta^t\left(\log\frac{R^N(s-)-1}{R^N(s-)}+\frac{1}{R^N(s-)}\right)dK_f^N(s)\\
&+\int_\beta^t\frac{1}{R^N(s-)}\,d\tilde K_a^N(s)+\int_\beta^t\frac{\rho-1}{R^N(s-)}\,d\tilde K_b^N(s)\\
&-\int_\beta^t\frac{\rho}{R^N(s-)}\,d\tilde K_d^N(s)-\int_\beta^t\frac{1}{R^N(s-)}\,d\tilde K_f^N(s)\\
&+\int_\beta^t\frac{(1-0.25\rho)X_1^N(s)+(\rho-1)2.5N^{-2/3}X_2^N(s)}{R^N(s)}\,ds\\
&-\int_\beta^t\frac{0.75N^{-2/3}X_2^N(s)Z_3^N(s)}{R^N(s)}\,ds.
\end{aligned}
$$

The second term on the right-hand-side is bounded by a constant times

$$
\begin{aligned}
\int_\beta^t&\left(\frac{K_a^N(s-)+1}{R^N(s-)}\right)^2\frac{1}{(K_a^N(s-)+1)^2}\,dK_a^N(s)\\
&\le\sup_{\beta\le s\le t}\left(\frac{K_a^N(s-)+1}{R^N(s-)}\right)^2\sum_{k=1}^\infty\frac{1}{k^2},
\end{aligned}
$$

and similarly for the third, fourth and fifth terms. By (3.6) and (3.8),

$$
\mathbf{1}_{\{\beta_{k_N}^N<\infty\}}\sup_{\beta_{k_N}^N\le s\le\tau_\varepsilon^N}\left(\frac{K_a^N(s-)+1}{R^N(s-)}\right)^2
$$

is stochastically bounded and similarly with $K_b^N$, $K_d^N$ and $K_f^N$. The sixth term is a martingale with quadratic variation

$$
\int_{\beta_{k_N}^N}^t\frac{1}{R^N(s-)^2}\,dK_a^N(s),
$$

and stochastic boundedness of the sequence of martingales follows from the stochastic boundedness of the quadratic variation and the boundedness of the jumps of the martingale. (See Lemma A.1.)

The last term satisfies

$$U^N(t) \equiv \int_\beta^t \frac{0.75N^{-2/3}X_2^N(s)Z_3^N(s)}{R^N(s)}\,ds$$
$$\leq 0.75N^{-2/3}\int_0^t Z_3^N(s)\,ds.$$

By the equations for $Z_3^N$,

$$E\left[\int_0^{\tau_\varepsilon^N} Z_3^N(s)\,ds\right] \leq 0.5E\left[\int_0^{\tau_\varepsilon^N} X_1^N(s)\,ds\right],$$

and adding the equations for $X_1^N$ and $X_2^N$ and taking expectations,

$$\begin{aligned}
&0.75E\left[\int_0^{\tau_\varepsilon^N} X_1^N(s)\,ds\right] + 1\\
&\quad = E[X_1^N(\tau_\varepsilon^N) + X_2^N(\tau_\varepsilon^N) + K_f^N(\tau_\varepsilon^N)]\\
&\quad \leq E[X_1^N(\tau_\varepsilon^N) + X_2^N(\tau_\varepsilon^N)] + 0.75\varepsilon E\left[\int_0^{\tau_\varepsilon^N} Z_3^N(s)\,ds\right]\\
&\quad \leq E[X_1^N(\tau_\varepsilon^N) + X_2^N(\tau_\varepsilon^N)] + 0.375\varepsilon E\left[\int_0^{\tau_\varepsilon^N} X_1^N(s)\,ds\right],
\end{aligned}$$

and hence,

$$\begin{aligned}
0.75(1-0.5\varepsilon)E\left[\int_0^{\tau_\varepsilon^N} X_1^N(s)\,ds\right] &\leq E[X_1^N(\tau_\varepsilon^N) + X_2^N(\tau_\varepsilon^N)]\\
&\leq N^{2/3}\varepsilon + \rho
\end{aligned} \tag{3.8}$$

and

$$E[U^N(\tau_\varepsilon^N)] \leq \frac{\varepsilon + N^{-2/3}\rho}{2-\varepsilon}.$$

Consequently,

$$\begin{aligned}
&\log R^N(\tau_\varepsilon^N)\\
&\quad = \log R^N(\beta_{k_N}^N) + O(1)\\
&\qquad + \int_{\beta_{k_N}^N}^{\tau_\varepsilon^N} \frac{(1-0.25\rho)X_1^N(s) + (\rho-1)2.5N^{-2/3}X_2^N(s)}{R^N(s)}\,ds\\
&\quad = \log R^N(\beta_{k_N}^N) + O(1) + \lambda^N(\tau_\varepsilon^N - \beta_{k_N}^N).
\end{aligned}$$

Note that

$$\begin{aligned}\limsup_{N\to\infty}\frac{\beta^N_{k_N}}{k_N N^{2/3}} &\le \limsup_{N\to\infty}\frac{1}{k_N N^{2/3}}\int_0^{\beta^N_{k_N}}(X_1^N(s)+X_2^N(s))\,ds\\ &\le \limsup_{N\to\infty}\frac{1}{2.5k_N}K_b^N(\beta_{k_N})\\ &\le \frac{1}{2.5}.\end{aligned}$$

Then assuming $k_N/\log N\to 0$,

$$\frac{\tau^N_\varepsilon-\beta^N_{k_N}}{N^{2/3}\log N}\approx\frac{\tau^N_\varepsilon}{N^{2/3}\log N},$$

and since $R^N(\beta^N_{k_N})\le(1+\rho^N)(1+k_N)$ and $N^{2/3}\lambda^N\to 7.5$, we have

$$\lambda^N N^{2/3}\frac{\tau^N_\varepsilon}{N^{2/3}\log N}\approx\frac{\log R^N(\tau^N_\varepsilon)}{\log N}\to\frac{2}{3},$$

giving the result. □

The computations in the proof of Theorem 3.3 also give the following lemma.

LEMMA 3.4. *Let* $0<\varepsilon_1<\varepsilon_2$. *Then* $N^{-2/3}(\tau^N_{\varepsilon_2}-\tau^N_{\varepsilon_1})$ *is stochastically bounded, and for* $\delta>0$,

$$\begin{aligned}\lim_{N\to\infty}P\biggl\{&\frac{2}{15}\log\frac{\varepsilon_2}{\varepsilon_1}-\delta\\ &\le N^{-2/3}(\tau^N_{\varepsilon_2}-\tau^N_{\varepsilon_1})\\ &\le\frac{2}{15}\log\frac{\varepsilon_2}{\varepsilon_1}+\int_{\tau^N_{\varepsilon_1}}^{\tau^N_{\varepsilon_2}}\frac{0.1N^{-2/3}X_2^N(s)Z_3^N(s)}{R^N(s)}\,ds+\delta\biggr\}=1.\end{aligned}$$

REMARK 3.5. In fact, we will see that $N^{-2/3}(\tau^N_{\varepsilon_2}-\tau^N_{\varepsilon_1})$ converges to a constant.

PROOF OF LEMMA 3.4. Since for $\tau^N_{\varepsilon_1}\le t\le\tau^N_{\varepsilon_2}$, $R_N(t)=O(N^{2/3})$, the integral expression for $\log R_N(t)$ implies

$$\log\frac{R_N(\tau^N_{\varepsilon_2})}{R_N(\tau_{\varepsilon_1})}=o(1)+\lambda^N(\tau^N_{\varepsilon_2}-\tau^N_{\varepsilon_1})-\int_{\tau^N_{\varepsilon_1}}^{\tau^N_{\varepsilon_2}}\frac{0.75N^{-2/3}X_2^N(s)Z_3^N(s)}{R^N(s)}\,ds.$$

The lemma follows from the fact that the last term is nonnegative and stochastically bounded. □

On $t \in [0, \tau_\varepsilon^N]$, $X_1^N(t)$ is dominated by the linear death process with immigration satisfying

$$\tilde{X}_1(t) = Y_b(2.5\varepsilon t) - Y_d\biggl(\int_0^t 0.25\tilde{X}_1(s)\,ds\biggr),$$

that is, an infinite server queue with Poisson arrivals of rate $2.5\varepsilon$ and exponential service times with rate 0.25. For $\beta > 0$, let $\gamma_\beta^N = \inf\{t : \tilde{X}_1(t) \geq N^\beta\}$. By Dynkin's formula, for each $t > 0$,

$$E[f(\tilde{X}_1(t \wedge \gamma_\beta^N))] = f(0) + E\biggl[\int_0^{\gamma_\beta^N \wedge t} Cf(\tilde{X}_1(s))\,ds\biggr],$$

where

$$Cf(k) = 2.5\varepsilon(f(k+1) - f(k)) + 0.25k(f(k-1) - f(k)).$$

For $0 < \delta < 1$, let $f(k) = (k!)^\delta$. Then $c_0 \equiv \sup_k Cf(k) < \infty$, and

$$P\{\gamma_\beta^N \leq t\} \leq \frac{1}{f(\lceil N^\beta \rceil)} E[f(\tilde{X}_1(t \wedge \gamma_\beta^N))] \leq \frac{1 + c_0 t}{f(\lceil N^\beta \rceil)}.$$

It follows that, for any $\alpha > 0$,

$$\lim_{N\to\infty} P\{\gamma_\beta^N \leq N^\alpha\} = 0.$$

Consequently, taking $0 < \beta < 2/3$ and $\alpha > 2/3$, we see that $N^{-2/3}X_1^N(\tau_\varepsilon^N) \to 0$ and, hence, $N^{-2/3}X_2^N(\tau_\varepsilon^N) \to \varepsilon$.

It is clear that the time scale for $Z_2^N$ is much slower than the time scale for $Z_1^N$ and $Z_3^N$. The fast time scale of $\{Z_1^N, Z_3^N\}$ enables us to "average out" their contributions to the evolution of the second component after time $\tau_\varepsilon^N$. We define $V_i^N(t) = Z_i(\tau_\varepsilon^N + N^{2/3}t)$.

Theorem 3.6 (Averaging and deterministic approximation).

(a) *Conditioning on $\tau_\varepsilon^N < \infty$, for each $\delta > 0$ and $t > 0$,*

$$\lim_{N\to\infty} P\Bigl\{\sup_{0\leq s\leq t} |V_2^N(s) - V_2(s)| \geq \delta\Bigr\} = 0,$$

*where $V_2$ is the solution of*

$$V_2(t) = \varepsilon + \int_0^t 7.5V_2(s)\,ds - \int_0^t 3.75V_2(s)^2\,ds. \tag{3.9}$$

(b) *Conditioning on $\tau_\varepsilon^N < \infty$, for each $t \geq 0$, $(V_1^N(t), V_3^N(t))$ converges in distribution to a pair $(V_1(t), V_3(t))$ with joint distribution $\mu_t^{13}$ satisfying*

$$\begin{aligned}\int \biggl[&2.5V_2(t)(g(z+1,y) - g(z,y)) \\ &+ 0.25z(g(z-1,y) - g(z,y)) + (z-2y)\frac{\partial g}{\partial y}(z,y)\biggr]\mu_t^{13}(dz,dy) = 0.\end{aligned} \tag{3.10}$$

*In particular,* $V_1(t)$ *has a Poisson distribution with parameter* $10V_2(t)$*, so*

$$E[V_1(t)] = \mathrm{Var}(V_1(t)) = 10V_2(t),$$
$$E[V_3(t)] = 5V_2(t), \qquad \mathrm{Var}(V_3(t)) = \tfrac{20}{9}V_2(t)$$

*and*

$$\mathrm{Cov}(V_1(t), V_3(t)) = \tfrac{40}{9}V_2(t).$$

REMARK 3.7. (a) Of course, the equation in Part (a) is just the classical logistic equation.

(b) For times $t_1 < t_2$, $(V_1^N(t_1), V_3^N(t_1))$ and $(V_1^N(t_2), V_3^N(t_2))$ are asymptotically independent.

PROOF OF THEOREM 3.6. On the event $\tau_\varepsilon^N < \infty$,

$$V_1^N(t) = Z_1^N(\tau_\varepsilon^N) + Y_b^\varepsilon\left(\int_0^t 2.5N^{2/3}V_2^N(s)\,ds\right) - Y_d^\varepsilon\left(\int_0^t 0.25N^{2/3}V_1^N(s)\,ds\right),$$

$$\begin{aligned} V_2^N(t) = {} & Z_2^N(\tau_\varepsilon^N) + N^{-2/3}Y_a^\varepsilon\left(\int_0^t N^{2/3}V_1^N(s)\,ds\right) \\ & - N^{-2/3}Y_b^\varepsilon\left(\int_0^t 2.5N^{2/3}V_2^N(s)\,ds\right) \\ & - N^{-2/3}Y_f^\varepsilon\left(N^{2/3}\int_0^t 0.75V_2^N(s)V_3^N(s)\,ds\right), \end{aligned}$$

$$\begin{aligned} V_3^N(t) = {} & Z_3^N(\tau_\varepsilon^N) + N^{-1}Y_c^\varepsilon\left(\int_0^t N^{5/3}V_1^N(s)\,ds\right) \\ & - N^{-1}Y_e^\varepsilon\left(\int_0^t 2N^{5/3}V_3^N(s)\,ds\right) \\ & - N^{-1}Y_f^\varepsilon\left(\int_0^t 0.75N^{2/3}V_2^N(s)V_3^N(s)\,ds\right), \end{aligned}$$

where $Y_a^\varepsilon$, $Y_b^\varepsilon$, and so on, are unit Poisson processes obtained from $Y_a$, $Y_b$, and so on. by taking increments over the appropriate intervals. For example,

$$Y_b^\varepsilon(u) = Y_b\left(\int_0^{\tau_\varepsilon^N} 2.5Z_2^N(s)\,ds + u\right) - Y_b\left(\int_0^{\tau_\varepsilon^N} 2.5Z_2^N(s)\,ds\right).$$

By the martingale properties of the Poisson processes,

$$\begin{aligned} & E[\rho^N N^{-2/3}V_1^N(t) + V_2^N(t)] \\ & \quad = E[\rho^N N^{-2/3}Z_1^N(\tau_\varepsilon^N) + Z_2^N(\tau_\varepsilon^N)] \\ & \qquad + \int_0^t E[(1 - 0.25\rho^N)V_1^N(s) + 2.5(\rho^N - 1)V_2^N(s)]\,ds \end{aligned}$$

$$
\begin{aligned}
&\quad - \int_0^t E[0.75 V_2^N(s) V_3^N(s)]\,ds \\
&\le E[\rho^N N^{-2/3} Z_1^N(\tau_\varepsilon^N) + Z_2^N(\tau_\varepsilon^N)] \\
&\quad + \lambda^N N^{2/3} \int_0^t E[\rho^N N^{-2/3} V_1^N(s) + V_2^N(s)]\,ds,
\end{aligned}
$$

so by Gronwall's lemma,

$$
\begin{aligned}
E[\rho^N N^{-2/3} V_1^N(t) + V_2^N(t)] &\le (\rho^N N^{-2/3} Z_1^N(\tau_\varepsilon^N) + Z_2^N(\tau_\varepsilon^N)) e^{N^{2/3}\lambda^N t} \\
&\le (\varepsilon + \rho^N N^{-2/3}) e^{N^{2/3}\lambda^N t},
\end{aligned}
$$

and the equation for $V_1^N$ implies

$$
\begin{aligned}
E[V_1^N(t)] &= E[Z_1^N(\tau_\varepsilon^N)] e^{-0.25 N^{2/3} t} + \int_0^t 2.5 N^{2/3} e^{-0.25 N^{2/3}(t-s)} E[V_2^N(s)]\,ds \\
&\le E[Z_1^N(\tau_\varepsilon^N)] e^{-0.25 N^{2/3} t} + 10(\varepsilon + \rho^N N^{-2/3}) e^{N^{2/3} t \lambda^N},
\end{aligned}
$$

and similarly for $V_3^N$,

$$
\begin{aligned}
E[V_3^N(t)] &\le E[Z_3^N(\tau_\varepsilon^N)] e^{-2N^{2/3} t} + \int_0^t N^{2/3} e^{-2N^{2/3}(t-s)} E[V_1^N(s)]\,ds \\
&\le E[Z_3^N(\tau_\varepsilon^N)] e^{-2N^{2/3} t} + E[Z_1^N(\tau_\varepsilon^N)] e^{-0.25 N^{2/3} t} \\
&\quad + 5(\varepsilon + \rho^N N^{-2/3}) e^{N^{2/3} t \lambda^N}.
\end{aligned}
$$

The law of large numbers for the Poisson process implies that $V_2^N$ is asymptotic to

$$
\widehat{V}_2^N(t) = Z_2^N(\tau_\varepsilon^N) + \int_0^t (V_1^N(s) - 2.5 V_2^N(s) - 0.75 V_2^N(s) V_3^N(s))\,ds,
$$

and since $Z_2^N(\tau_\varepsilon^N) \to \varepsilon$, as in the proof of Lemma A.4, we can verify the relative compactness of $\{\widehat{V}_2^N\}$ if we can verify the stochastic boundedness of

$$
\int_0^t (V_1^N(s) - 2.5 V_2^N(s) - 0.75 V_2^N(s) V_3^N(s))^2\,ds,
$$

which in turn will follow from the stochastic boundedness of

$$
\int_0^t V_1^N(s)^k\,ds, \qquad \int_0^t V_2^N(s)^k\,ds, \qquad \int_0^t V_3^N(s)^k\,ds, \tag{3.11}
$$

for appropriate $k$.

Let $\gamma^N = \inf\{t : Z_2^N(t) > 4\}$, and in the equation for $Z_1^N$, replace $Y_b(\int_0^t 2.5 \times Z_2^N(s)\,ds)$ by $Y_b(\int_0^{t\wedge\gamma^N} 2.5 Z_2^N(s)\,ds)$. Note that $\gamma^N > \tau_\varepsilon^N$, and if we can verify the stochastic boundedness of (3.11) for the modified system and show that

$\gamma^N \Rightarrow \infty$, we will have the stochastic boundedness for the original system. Note that

$$\int_0^t V_1^N(s)^k\,ds = N^{-2/3}\int_{\tau_\varepsilon^N}^{\tau_\varepsilon^N+N^{2/3}t} Z_1^N(u)^k\,du. \tag{3.12}$$

Taking $\varepsilon_0 = \frac{\varepsilon}{2}$ and $t_0 < \frac{2}{15}\log 2$, it follows from Lemma 3.4 that the sequence in (3.12) is stochastically bounded for each $t$ if and only if

$$N^{-2/3}\int_{\tau_{\varepsilon_0}^N+N^{2/3}t_0}^{\tau_{\varepsilon_0}^N+N^{2/3}t_1} Z_1^N(u)^k\,du = \int_{t_0}^{t_1} Z_1^N(\tau_{\varepsilon_0}^N + N^{2/3}s)^k\,ds$$

is stochastically bounded for each $t_1$. By Lemma A.3,

$$\begin{aligned} &E[Z_1^N(\tau_{\varepsilon_0}^N + N^{2/3}s)^k|\mathcal{F}_{\tau_{\varepsilon_0}^N}^N] \\ &\quad\le (Z_1^N(\tau_{\varepsilon_0}^N)^k \vee 1)e^{-0.24N^{2/3}s} + C_k 11(1-e^{-0.25N^{2/3}s}), \end{aligned}\tag{3.13}$$

where the 11 comes from the fact that $2.5Z_2^N(u)\mathbf{1}_{\{u<\gamma^N\}} + 1 \le 11$. Since $Z_1^N(\tau_{\varepsilon_0}^N) \le \frac{N^{2/3}\varepsilon_0}{\rho_N} + 1$, the first term goes to zero, and the stochastic boundedness follows. Stochastic boundedness for

$$\int_{t_0}^{t_1} Z_i^N(\tau_{\varepsilon_0}^N + N^{2/3}s)^k\,ds,$$

$i = 2, 3$, follows by a similar argument, using the fact that $Z_2^N(\tau_{\varepsilon_0}^N) \le \varepsilon_0 + N^{-2/3}$ and, by (3.8),

$$\begin{aligned} E[Z_3^N(\tau_{\varepsilon_0}^N)] &\le E\left[\int_0^{\tau_{\varepsilon_0}^N} X_1^N(s)\,ds\right] \\ &\le \frac{N^{2/3}\varepsilon_0 + \rho_N}{0.75(1-0.5\varepsilon_0)}, \end{aligned}$$

and applying (3.13) to bound the second term on the right-hand-side of (A.5).

As $N \to \infty$, dividing the equations for $V_1^N$ and $V_3^N$ by $N^{2/3}$ shows that

$$\int_0^t V_1^N(s)\,ds - 10\int_0^t V_2^N(s)\,ds \to 0,$$

$$\int_0^t V_3^N(s)\,ds - 5\int_0^t V_2^N(s)\,ds \to 0.$$

The assertion for $V_3^N$ and the fact that $V_2^N$ is asymptotically regular (e.g., one can prove that $\lim_{h\to 0}\limsup_{N\to\infty} E[\sup_{t\le T}\sup_{s\le h}|V_2^N(t+s) - V_2^N(t)|] = 0$) implies

$$\int_0^t V_2^N(s)V_3^N(s)\,ds - 5\int_0^t V_2^N(s)^2\,ds \to 0.$$

It follows that $V_2^N$ converges to the solution of (3.9). It should now be clear why we shifted the initial time to $\tau_\varepsilon^N$.

$V_1^N$ and $V_3^N$ fluctuate rapidly and locally in time. $V_1^N$ behaves like a simple birth and death process with $V_2^N$ entering as a parameter, and $V_3^N$ can be approximated by an ordinary differential equation driven by $V_1^N$, that is,

$$V_3^N(a+N^{-2/3}r) \approx V_3^N(a)e^{-2r} + \int_0^r e^{-2(r-s)} V_1^N(a+N^{-2/3}s)\,ds.$$

To be specific, let

$$\begin{aligned} B_s^N g(z,y) = {} & 2.5V_2^N(s)(g(z+1,y)-g(z,y)) \\ & + 0.25z(g(z-1,y)-g(z,y)) \\ & + zN(g(z,y+N^{-1})-g(z,y)) \\ & + 2yN(g(z,y-N^{-1})-g(z,y)). \end{aligned}$$

Then

$$\begin{aligned} M_g^N(t) \equiv {} & g(V_1^N(t),V_3^N(t)) - g(Z_1(\tau_\varepsilon^N),Z_3(\tau_\varepsilon^N)) \\ & - N^{2/3}\int_0^t B_s^N g(V_1^N(s),V_3^N(s))\,ds \end{aligned}$$

is a martingale, and defining an occupation measure $\Gamma^N$ by

$$\Gamma^N(C\times D\times[0,t]) = \int_0^t \mathbf{1}_C(V_1^N(s))\mathbf{1}_D(V_3^N(s))\,ds,$$

$$\begin{aligned} M_g^N(t) = {} & g(V_1^N(t),V_3^N(t)) - g(Z_1(\tau_\varepsilon^N),Z_3(\tau_\varepsilon^N)) \\ & - N^{2/3}\int_{\mathbb{Z}^+\times\mathbb{R}^+\times[0,t]} B_s^N g(z,y)\Gamma^N(dz\times dy\times ds). \end{aligned} \tag{3.14}$$

Let $\mathcal{L}_m \equiv \mathcal{L}_m(\mathbb{Z}^+\times\mathbb{R}^+)$ denote the space of measures $\nu$ on $\mathbb{Z}^+\times\mathbb{R}^+\times[0,\infty)$ such that $\nu(\mathbb{Z}^+\times\mathbb{R}^+\times[0,t]) = t$, topologized so that convergence is weak convergence on $\mathbb{Z}^+\times\mathbb{R}^+\times[0,t]$ for each $t>0$. It is easy to verify that the sequence $(V_2^N,\Gamma^N)$ is relatively compact in $D_{\mathbb{R}^+}([0,\infty))\times\mathcal{L}_m$, where $D_{\mathbb{R}^+}([0,\infty))$ is the space of cadlag $\mathbb{R}^+$-valued functions. Let $(V_2,\Gamma)$ be a limit point of the sequence.

Define

$$\begin{aligned} B_v g(z,y) = {} & 2.5v(g(z+1,y)-g(z,y)) + 0.25z(g(z-1,y)-g(z,y)) \\ & + (z-2y)\frac{\partial g}{\partial y}(z,y). \end{aligned}$$

Noting that

$$\lim_{N\to\infty}(B_s^N g(z,y) - B_{V_2^N(s)}g(z,y)) = 0$$

and dividing (3.14) by $N^{2/3}$ and letting $N \to \infty$, we have

$$\int_{\mathbb{Z}^+ \times \mathbb{R}^+ \times [0,t]} B_{V_2(s)} g(z,y) \Gamma(dz \times dy \times ds) = 0. \tag{3.15}$$

Differentiating (3.15) gives (3.10) at least for almost every $t$. (See [8] for more details.)

From (3.10), we can easily obtain all the moments of the limiting joint distribution. Let $(Z_s, Y_s)$ be a random vector with the law $\mu_s^{13}$. Note that $Z_s$ is just a Poisson random variable with expectation $10V_2(s)$. Hence, the marginal moments of $Z_s$ are immediate. With $g(z,y) = y$ in (3.10), we get

$$\int [z - 2y] \mu_s^{13}(dz, dy) = 0$$

and, hence,

$$E[Y_s] = \tfrac{1}{2} E[Z_s] = 5V_2(s).$$

With $g(z,y) = zy$, we get

$$\int [z(z-2y) + \tfrac{5}{2} y V_2(s) - \tfrac{1}{4} yz] \mu_s^{13}(dz, dy) = 0,$$

$$E[Z_s Y_s] = \tfrac{4}{9} E[Z_s^2] + \tfrac{10}{9} V_2(s) E[Y_s] = \tfrac{40}{9} V_2(s) + 50 V_2(s)^2.$$

With $g(z,y) = y^2$,

$$\int [2y(z-2y)] \mu_s^{13}(dz, dy) = 0$$

and

$$E[Y_s^2] = \tfrac{1}{2} E[Z_s Y_s].$$

In general, taking $g(z,y) = z^n y^m$, one gets the recursive equation

$$\begin{aligned} & mE[Z_s^{n+1} Y_s^{m-1}] - 2mE[Z_s^n Y_s^m] + \tfrac{5}{2} V_2(s) \sum_1^{n-1} \binom{n}{k} E[Z_s^{n-k} Y_s^m] \\ & \quad + \tfrac{1}{4} \sum_1^{n-1} \binom{n}{k} (-1)^k E[Z_s^{n-k+1} Y_s^m] = 0, \end{aligned} \tag{3.16}$$

from which one can iteratively compute all the joint moments of $(Z_s, Y_s)$. Note that, in order to compute $E[Z_s^n Y_s^m]$, one first has to compute

$$E[Z_s^{n+m}], \qquad E[Z_s^{n+m-1} Y_s], \qquad E[Z_s^{n+m-2} Y_s^2], \ldots, E[Z_s^{n+1} Y_s^{m-1}],$$

as well as all of $\{E[Z_s^i Y_s^j] : 0 < i + j < n + m\}$. $\square$

# APPENDIX

## A.1. Estimates.

LEMMA A.1. *Suppose $M$ is a martingale with quadratic variation $[M]$ and $Z = \sup_t |M(t) - M(t-)|$. Then*

$$P\Big\{\sup_{s\le t}|M(s)| \ge c\Big\} \le \frac{d + E[Z^2]}{c^2} + P\{[M]_t > d\}.$$

PROOF. Let $\tau = \inf\{s : [M]_s > d\}$ and note that $[M]_{\tau\wedge t} \le d + Z^2$. Then

$$P\Big\{\sup_{s\le t}|M(s)| \ge c\Big\} \le P\Big\{\sup_{s\le t}|M(\tau\wedge s)| \ge c\Big\} + P\{\tau \le t\},$$

and the result follows by Doob's inequality. □

We need the following inequality.

LEMMA A.2. *Let $z(t) \ge 0$, $t \ge 0$, $F : [0,\infty) \to [0,\infty)$ be nondecreasing, and $b > 0$. Suppose that, for $0 \le r < t$*

$$z(t) - z(r) \le \int_r^t (F(z(s)) - bz(s))\,ds. \tag{A.1}$$

*If $z(0) \le z^*(0)$ and $z^*$ satisfies*

$$z^*(t) = z^*(0) + \int_0^t (F(z^*(s)) - bz^*(s))\,ds,$$

*then $z(t) \le z^*(t)$, $t \ge 0$. If, in addition, $z^*(0) \ge 1$, $c \ge b$, and for some $k > 1$, $F(u) \le cu^{(k-1)/k}$, $u \ge 1$, then*

$$\begin{aligned} z(t) &\le \Big(z^*(0)^{1/k} e^{-(b/k)t} + \frac{c}{b}(1 - e^{-(b/k)t})\Big)^k \\ &\le z^*(0)e^{-(b/k)t} + \frac{c^k}{b^k}(1 - e^{-(b/k)t}). \end{aligned} \tag{A.2}$$

PROOF. Let $z_0 \equiv z$, and let $z_1$ satisfy

$$z_1(t) = z(0) + \int_0^t (F(z(s)) - bz(s))\,ds - \int_0^t b(z_1(s) - z(s)) \vee 0\,ds. \tag{A.3}$$

Let $\Gamma = \{r : z_1(r) \ge z(r)\}$. For $r \in \Gamma$, (A.1) and $t > r$,

$$z(t) \le z_1(r) + \int_r^t (F(z(s)) - bz(s))\,ds. \tag{A.4}$$

Let $t^* = \sup\{s < t : z_1(s) \geq z(s)\}$. The continuity of $z_1$ and (A.4) imply that $z(t^*) \leq z_1(t^*)$. If $t > t^*$, then $z(s) > z_1(s)$, for $t^* < s < t$, and the last term in (A.3) is zero, so

$$z(t) \leq z_1(t^*) + \int_{t^*}^t (F(z(s)) - bz(s))\,ds = z_1(t).$$

Consequently, $z_1(t) \geq z_0(t)$, for all $t \geq 0$, and

$$z_1(t) = z(0) + \int_0^t (F(z_0(s)) - bz_1(s))\,ds.$$

For $m > 0$, define $z_{m+1}$ recursively by requiring

$$z_{m+1}(t) = z(0) + \int_0^t (F(z_m(s)) - bz_{m+1}(s))\,ds.$$

Then $z_{m+1} \geq z_m$ and $z_m$ converges to the solution of

$$z^*(t) = z(0) + \int_0^t (F(z^*(s)) - bz^*(s))\,ds.$$

It follows that $z^* \geq z$.

Suppose $F(u) \leq cu^{(k-1)/k}$, $u \geq 1$, $c > b$, $z^*(0) \geq 1$, and $z^{**}(t)$ satisfies

$$z^{**}(t) = z^*(0) + \int_0^t (cz^{**}(s)^{(k-1)/k} - bz^{**}(s))\,ds.$$

Note that $z^{**}(t) \geq 1$, $t \geq 0$, so $F(z^{**}(t)) \leq cz^{**}(t)^{(k-1)/k}$, $t \geq 0$, and $z^*(t) \leq z^{**}(t)$, $t \geq 0$. Setting $z^{**}(t) = y(t)^k$,

$$y(t) = z^*(0)^{1/k} + \frac{c}{k}t - \int_0^t \frac{b}{k} y(s)\,ds$$

and

$$y(t) = z^*(0)^{1/k} e^{-(b/k)t} + \frac{c}{b}(1 - e^{-(b/k)t}). \qquad \square$$

LEMMA A.3. *Let $\alpha, \beta \geq 0$. Suppose $Z \geq 0$ is adapted to a filtration $\{\mathcal{F}_t\}$ and*

$$Z(t) = N^{-\beta}(K_0 + K_1(t) - K_2(t) - K_3(t)),$$

*where $K_0$ is an integer-valued random variable and $K_1$, $K_2$ and $K_3$ are counting processes with no simultaneous jumps and $\{\mathcal{F}_t\}$-intensities $N^\alpha U(t)$, $\mu N^\alpha Z(t)$ and $\lambda_N(t)$, respectively. Then, there exists $C_k > 0$ such that*

$$\begin{aligned} E[Z(t)^k | Z(0)] &\leq (Z(0)^k \vee 1) e^{-\mu N^{\alpha-\beta} t} \\ &\quad + C_k \sup_{s \leq t} E[(U(s)+1)^k](1 - e^{-\mu N^{\alpha-\beta} t}). \end{aligned} \tag{A.5}$$

PROOF. Let $Y$ be a unit Poisson process that is independent of the $K_i$, and let $\widehat{Z}$ satisfy

$$\widehat{Z}(t) = N^{-\beta}\left(K_0 + K_1(t) - K_2(t) - Y\left(\int_0^t \mu N^{\alpha}(\widehat{Z}(s) - Z(s))\,ds\right)\right).$$

Then $Z(t) \le \widehat{Z}(t)$ and

$$\widehat{K_2}(t) = K_2(t) + Y\left(\int_0^t \mu N^{\alpha}(\widehat{Z}(s) - Z(s))\,ds\right)$$

has intensity $\mu N^{\alpha}\widehat{Z}(t)$. We have

$$\begin{aligned}\widehat{Z}(t)^k = Z(0)^k &+ \int_0^t ((\widehat{Z}(s-) + N^{-\beta})^k - \widehat{Z}(s-)^k)\,dK_1(s)\\ &+ \int_0^t ((\widehat{Z}(s-) - N^{-\beta})^k - \widehat{Z}(s-)^k)\,d\widehat{K_2}(s)\end{aligned}$$

and

$$\begin{aligned}&E[\widehat{Z}(t)^k | Z(0)]\\ &\quad = Z(0)^k + \int_0^t N^{\alpha}\sum_{l=0}^{k-1}\binom{k}{l} E[U(s)\widehat{Z}(s)^l N^{-\beta(k-l)}\\ &\qquad\qquad + \mu(-1)^{k-l}\widehat{Z}(s)^{l+1}N^{-\beta(k-l)}|Z(0)]\,ds\\ &\quad = Z(0)^k + N^{\alpha-\beta}\int_0^t \sum_{l=0}^{k-1} E[H_{N,k,l}(U(s), \widehat{Z}(s)|Z(0))]\,ds\\ &\qquad - \int_0^t \mu N^{\alpha-\beta} k E[\widehat{Z}(s)^k|Z(0)]\,ds,\end{aligned}$$

where

$$\begin{aligned}&H_{N,k,l}(U(s), \widehat{Z}(s))\\ &\quad = \left(\binom{k}{l} U(s) - \mu(-1)^{k-l}\binom{k}{l-1} N^{-\beta}\right) N^{-\beta(k-l-1)}\widehat{Z}(s)^l.\end{aligned}$$

Applying the Hölder inequality, there exists $a_{k,l} > 0$ not depending on $N$ such that

$$\begin{aligned}&E[H_{N,k,l}(U(s), \widehat{Z}(s))|Z(0)]\\ &\qquad \le a_{k,l} E[(U(s)+1)^k|Z(0)]^{1/k} E[\widehat{Z}(s)^k|Z(0)]^{l/k}.\end{aligned}$$

It follows that

$$E[\widehat{Z}(t)^k | Z(0)]$$

$$\leq Z(0)^k + N^{\alpha-\beta}\int_0^t E[(U(s)+1)^k|Z(0)]^{1/k} G(E[\widehat{Z}(s)^k|Z(0)])\,ds$$
$$- \int_0^t \mu N^{\alpha-\beta} k E[\widehat{Z}(s)^k|Z(0)]\,ds,$$

where $G(u)$ is a polynomial of degree $k-1$ in $u^{1/k}$. It follows that, for $u \geq 1$, there is a constant $a_k$ such that $G(u) \leq a_k u^{(k-1)/k}$. Applying (A.2) with $b = \mu N^{\alpha-\beta} k$ and appropriate choice of $c$ gives (A.5). □

LEMMA A.4. *Suppose that $\beta > 0$ and that $Z^N(t) = N^{-\beta} K^N(t)$, where $K^N$ is a counting process with intensity $N^\beta \lambda_N(t)$. Suppose that $\{\int_0^t \lambda_N(s)^2\,ds\}$ is stochastically bounded for each $t > 0$. Then $\{Z^N\}$ is relatively compact as a sequence of processes in the sense of convergence in distribution in $D_{\mathbb{R}}[0,\infty)$, and every limit point has continuous sample paths.*

PROOF. Since $K^N$ can be represented as $Y(N^\beta \int_0^t \lambda_N(s)\,ds)$ for a unit Poisson process $Y$, by the law of large numbers, it is enough to verify the relative compactness of the sequence $\Lambda^N(t) = \int_0^t \lambda_N(s)\,ds$. But for $t < t+h \leq T$,

$$|\Lambda^N(t+h) - \Lambda^N(t)| \leq \sqrt{h}\sqrt{\int_0^T \lambda_N(s)^2\,ds},$$

which gives a uniform equicontinuity condition implying the relative compactness of $\{\Lambda^N\}$ (see, e.g., Theorem 3.7.2 of [4]). □

**A.2. Determining the scaling exponents.** The scalings employed for the examples in Sections 2 and 3 were determined in part by examining the published simulations. In particular, these simulations suggested the relationships among the $\alpha_k$. That approach to the choice of the scalings, however, is still more art than science and leaves open the question of whether slightly different, but equally reasonable scalings would produce significantly different limiting approximations. In this section we reconsider the model of Section 3 and give a more systematic identification of the scaling.

Recall that the basic model satisfies

$$X_1(t) = X_1(0) + Y_b\left(\int_0^t \kappa_2 X_2(s)\,ds\right) - Y_d\left(\int_0^t \kappa_4 X_1(s)\,ds\right),$$
$$X_2(t) = X_2(0) + Y_a\left(\int_0^t \kappa_1 X_1(s)\,ds\right) - Y_b\left(\int_0^t \kappa_2 X_2(s)\,ds\right)$$
$$- Y_f\left(\int_0^t \kappa_6 X_2(s) X_3(s)\,ds\right),$$
$$X_3(t) = X_3(0) + Y_c\left(\int_0^t \kappa_3 X_1(s)\,ds\right) - Y_e\left(\int_0^t \kappa_5 X_3(s)\,ds\right)$$

$$- Y_f\left(\int_0^t \kappa_6 X_2(s)X_3(s)\,ds\right).$$

With reference to Section 1.2, we consider a general scaling $Z_i^N(t) = N^{-\alpha_i}X_i(t)$, and replace $\kappa_k$ by $\lambda_k N^{\beta_k}$. Once the $\beta_k$ are selected, the $\lambda_k$ are determined by setting

$$\lambda_k = \kappa_k N_0^{-\beta_k}$$

for the rate constants $\kappa_k$ given in Section 3 and some appropriate $N_0$.

The normalized system becomes

$$\begin{aligned}
Z_1^N(t) &= Z_1^N(0) + N^{-\alpha_1}Y_b\left(\int_0^t \kappa_2 N^{\alpha_2} Z_2^N(s)\,ds\right)\\
&\quad - N^{-\alpha_1}Y_d\left(\int_0^t \kappa_4 N^{\alpha_1} Z_1^N(s)\,ds\right),\\
Z_2^N(t) &= Z_2^N(0) + N^{-\alpha_2}Y_a\left(\int_0^t \kappa_1 N^{\alpha_1} Z_1^N(s)\,ds\right)\\
&\quad - N^{-\alpha_2}Y_b\left(\int_0^t \kappa_2 N^{\alpha_2} Z_2^N(s)\,ds\right)\\
&\quad - N^{-\alpha_2}Y_f\left(\int_0^t \kappa_6 N^{\alpha_2+\alpha_3} Z_2^N(s)Z_3^N(s)\,ds\right),\\
Z_3^N(t) &= Z_3^N(0) + N^{-\alpha_3}Y_c\left(\int_0^t \kappa_3 N^{\alpha_1} Z_1^N(s)\,ds\right)\\
&\quad - N^{-\alpha_3}Y_e\left(\int_0^t \kappa_5 N^{\alpha_3} Z_3^N(s)\,ds\right)\\
&\quad - N^{-\alpha_3}Y_f\left(\int_0^t \kappa_6 N^{\alpha_2+\alpha_3} Z_2^N(s)Z_3^N(s)\,ds\right).
\end{aligned}$$

Setting $V_k^N(t) = Z_k^N(N^\gamma t)$ and replacing $\kappa_i$ by $\lambda_i N^{\beta_i}$,

$$\begin{aligned}
V_1^N(t) &= V_1^N(0) + N^{-\alpha_1}Y_b\left(\int_0^t \lambda_2 N^{\gamma+\beta_2} N^{\alpha_2} V_2^N(s)\,ds\right)\\
&\quad - N^{-\alpha_1}Y_d\left(\int_0^t \lambda_4 N^{\gamma+\beta_4} N^{\alpha_1} V_1^N(s)\,ds\right),\\
V_2^N(t) &= V_2^N(0) + N^{-\alpha_2}Y_a\left(\int_0^t \lambda_1 N^{\gamma+\beta_1} N^{\alpha_1} V_1^N(s)\,ds\right)\\
&\quad - N^{-\alpha_2}Y_b\left(\int_0^t \lambda_2 N^{\gamma+\beta_2} N^{\alpha_2} V_2^N(s)\,ds\right)\\
&\quad - N^{-\alpha_2}Y_f\left(\int_0^t \lambda_6 N^{\gamma+\beta_6} N^{\alpha_2+\alpha_3} V_2^N(s)V_3^N(s)\,ds\right),
\end{aligned}$$

$$
\begin{aligned}
V_3^N(t) = V_3^N(0) &+ N^{-\alpha_3} Y_c\left(\int_0^t \lambda_3 N^{\gamma+\beta_3} N^{\alpha_1} V_1^N(s)\,ds\right)\\
&- N^{-\alpha_3} Y_e\left(\int_0^t \lambda_5 N^{\gamma+\beta_5} N^{\alpha_3} V_3^N(s)\,ds\right)\\
&- N^{-\alpha_3} Y_f\left(\int_0^t \lambda_6 N^{\gamma+\beta_6} N^{\alpha_2+\alpha_3} V_2^N(s) V_3^N(s)\,ds\right).
\end{aligned}
$$

We assume that $(V_1^N(0), V_2^N(0), V_3^N(0)) \to (V_1(0), V_2(0), V_3(0))$.

The question is how to determine, in a systematic way, what the exponents $\alpha_i$, $\beta_k$ and $\gamma$ should be. There are several conditions that help this determination. First, we want the scaling to ensure that $V_i^N(t) = O(1)$. This requirement can be met either by ensuring that the individual terms on the right are $O(1)$ or by ensuring that terms cancel. Second, it is natural to assume that the $\beta_k$ have the same order as the $\kappa_k$, that is, we should have

$$
\beta_6 \le \beta_2 \le \beta_4 \le \beta_1 \le \beta_5 \le \beta_3. \tag{A.6}
$$

As is clearly reasonable, we assume that $\beta_1 = 0$. This last assumption is not really a restriction, since if $\beta_1 \ne 0$, we can add $\beta_1$ to $\gamma$ and substract $\beta_1$ from each of the $\beta_k$.

Finally, comparing the $\kappa_k$, it is also natural to assume that $\beta_3 > \beta_5$ and $\beta_6 < \beta_2$. (We will see that the second of these assumptions is actually implied by other considerations.) For the scaling used in Section 3, $\beta_1 = \beta_4 = \beta_5 = 0$, $\beta_2 = -2/3$, $\beta_3 = 1$ and $\beta_6 = -5/3$.

Suppose, as is the case in Section 3, we also require that the scaling makes each of the terms in the equation for $V_2^N$ to be $O(1)$. In particular, we look for a scaling in which the nonlinear behavior is preserved. Then we must have

$$
\begin{aligned}
\alpha_2 &= \gamma + \alpha_1,\\
\alpha_2 &= \gamma + \beta_2 + \alpha_2,\\
\alpha_2 &= \gamma + \beta_6 + \alpha_2 + \alpha_3.
\end{aligned}
$$

In addition, for $V_1^N$ to be $O(1)$ without being asymptotically negligible, we must have

$$
\alpha_2 = \gamma + \beta_2 + \alpha_2 = \gamma + \beta_4 + \alpha_1.
$$

Similarly, for $V_3^N$, we must have

$$
\begin{aligned}
\gamma + \beta_3 + \alpha_1 &\ge \gamma + \beta_5 + \alpha_3,\\
\gamma + \beta_3 + \alpha_1 &\ge \gamma + \beta_6 + \alpha_2 + \alpha_3 = \alpha_2,
\end{aligned}
$$

with equality holding for at least one of the inequalities. Since we are assuming that $\beta_3 > \beta_5 \ge 0$, we must have

$$
\gamma + \beta_3 + \alpha_1 = \gamma + \beta_5 + \alpha_3
$$

and, hence, $\alpha_3 > \alpha_1$.

The above assumptions imply

$$\alpha_2 - \alpha_1 = \beta_4 - \beta_2 = \gamma = -\beta_2 \geq 0,$$

so $\alpha_2 \geq \alpha_1$, $\beta_4 = 0$, and

$$\alpha_3 = -\gamma - \beta_6 = \beta_2 - \beta_6.$$

These restrictions leave three cases of interest: $\alpha_1 = \alpha_2 > 0$ and $\beta_5 = 0$, $\alpha_1 = \alpha_2 > 0$ and $\beta_5 > 0$, and $\alpha_2 > \alpha_1 \geq 0$.

If $\alpha_1 = \alpha_2 > 0$, then $\beta_2 = \beta_1 = \beta_4 = 0$ and $\alpha_3 = -\beta_6$. If, in addition, $\beta_5 = 0$, then $\alpha_3 = \gamma + \beta_3 + \alpha_1$ and as $N \to \infty$, the system converges to the solution of

$$V_1(t) = V_1(0) + \int_0^t (\lambda_2 V_2(s) - \lambda_4 V_1(s))\, ds,$$

$$V_2(t) = V_2(0) + \int_0^t (\lambda_1 V_1(s) - \lambda_2 V_2(s) - \lambda_6 V_2(s) V_3(s))\, ds,$$

$$V_3(t) = V_3(0) + \int_0^t (\lambda_3 V_1(s) - \lambda_5 V_3(s))\, ds.$$

If $\alpha_1 = \alpha_2 > 0$ and $\beta_5 > 0$, then

$$\lim_{N\to\infty} \left( \int_0^t (\lambda_3 V_1^N(s) - \lambda_5 V_3^N(s))\, ds \right) = 0, \tag{A.7}$$

and $(V_1^N, V_2^N)$ converges to a solution of

$$V_1(t) = V_1(0) + \int_0^t (\lambda_2 V_2(s) - \lambda_4 V_1(s))\, ds$$

$$V_2(t) = V_2(0) + \int_0^t \left( \lambda_1 V_1(s) - \lambda_2 V_2(s) - \frac{\lambda_6 \lambda_3}{\lambda_5} V_2(s) V_1(s) \right) ds.$$

If $\alpha_2 > \alpha_1$, then $\beta_2 < 0$ and $\gamma > 0$ and, consequently, $\gamma + \beta_4 > 0$ and $\gamma + \beta_5 > 0$. It follows that

$$\lim_{N\to\infty} \left( \int_0^t (\lambda_4 V_1^N(s) - \lambda_2 V_2^N(s))\, ds \right) = 0 \tag{A.8}$$

and (A.7) hold. Then, as in the scaling in Section 3, $V_2^N$ converges to the solution of

$$V_2(t) = V_2(0) + \int_0^t \left( \left( \frac{\lambda_1 \lambda_2}{\lambda_4} - \lambda_2 \right) V_2(s) - \frac{\lambda_6 \lambda_3 \lambda_2}{\lambda_5 \lambda_4} V_2(s)^2 \right) ds. \tag{A.9}$$

Define $\lambda_k = \kappa_k N_0^{-\beta_k}$ for some $N_0$. Then

$$\frac{\lambda_1 \lambda_2}{\lambda_4} - \lambda_2 = \left( \frac{\kappa_1 \kappa_2}{\kappa_4} - \kappa_2 \right) N_0^{\gamma}$$

and

$$\begin{aligned}\frac{\lambda_6\lambda_3\lambda_2}{\lambda_5\lambda_4} &= \frac{\kappa_6\kappa_3\kappa_2}{\kappa_5\kappa_4}N^{-\beta_6-\beta_3-\beta_2+\beta_5}\\ &= \frac{\kappa_6\kappa_3\kappa_2}{\kappa_5\kappa_4}N_0^{\alpha_2+\gamma}.\end{aligned}$$

Recalling that $V_2^{N_0}(t) = N_0^{-\alpha_2}X_2(N^\gamma t)$, the convergence suggests approximating $X_2(t)$ by $\widehat{V}_2(t) = N_0^{\alpha_2}V_2(N_0^{-\gamma}t)$. But if $V_2$ satisfies (A.9), then $\widehat{V}_2$ satisfies

$$\widehat{V}_2(t) = \widehat{V}_2(0) + \int_0^t \left(\left(\frac{\kappa_1\kappa_2}{\kappa_4} - \kappa_2\right)\widehat{V}_2(s) - \frac{\kappa_6\kappa_3\kappa_2}{\kappa_5\kappa_4}\widehat{V}_2(s)^2\right)ds, \tag{A.10}$$

so the approximation does not depend on the choice of the scaling parameters beyond the restrictions identified above and the assumption that $\alpha_2 > \alpha_1$.

The behavior of the $V_1^N$ and $V_3^N$ depends primarily on whether $\alpha_1 > 0$ or $\alpha_1 = 0$. If $\alpha_1 > 0$, then (A.7) and (A.8) can be strengthened to

$$\lim_{N\to\infty}\sup_{\varepsilon\le s\le t}(|\lambda_3 V_1^N(s) - \lambda_5 V_3^N(s)| + |\lambda_4 V_1^N(s) - \lambda_2 V_2^N(s)|) = 0$$

for each $0 < \varepsilon < t$.

If $\alpha_1 = 0$, then the behavior of $V_1^N$ is essentially the same as in Section 3. If, in addition, $\beta_5 = 0$, then the joint behavior of $V_1^N$ and $V_3^N$ is essentially the same as in Section 3. If $\alpha_1 = 0$ and $\beta_5 > 0$, then

$$\lim_{N\to\infty}\int_0^t |\lambda_3 V_1^N(s) - \lambda_5 V_3^N(s)|\,ds = 0, \tag{A.11}$$

and for each $t > 0$, $(V_1^N(t), V_3^N(t)) \Rightarrow (V_1(t), \frac{\lambda_3}{\lambda_5}V_1(t))$, where $V_1(t)$ is Poisson distributed with parameter $\frac{\lambda_2 V_2(t)}{\lambda_4}$.

**Acknowledgments.** The authors would like to thank the Institute of Mathematics and its Applications for their support and hospitality during the year on "Probability and Statistics in Complex Systems."

K. Ball
IDA Center for Communications Research
4320 Westerra Court
San Diego, California 92122
USA
E-mail: ball@ccrwest.org

T. G. Kurtz
Department of Mathematics
University of Wisconsin
480 Lincoln Drive
Madison, Wisconsin 53706
USA
E-mail: kurtz@math.wisc.edu

L. Popovic
Institute of Mathematics and its Applications
University of Minnesota
207 Church St SE
Minneapolis, Minnesota 55455
USA
E-mail: lea@ima.umn.edu

G. Rempala
Department of Mathematics
University of Louisville
Louisville, Kentucky 40292
USA
E-mail: grempala@louisville.edu